\documentclass[12pt]{amsart}
\usepackage{amsfonts,amsmath,amssymb,amscd}
\usepackage{amsthm}
\usepackage{enumerate}
\usepackage{bbm}

\usepackage[colorlinks = true,
            linkcolor = red,
            urlcolor  = blue,
            citecolor = blue,
            anchorcolor = blue]{hyperref}

\def\cF{\mathcal F}
\def\cK{\mathcal K}
\def\cS{\mathcal S}
\def\cL{{\mathcal L}} 
\def\cR{\mathcal R}

\def\Ga{\Gamma}
\def\bd{\partial} 
\def\Leb{{\rm Leb}}
\def\and{{\rm \text{ and }}}
\def\dist{{\rm \text{ dist }}}
\def\HD{{\text {\rm HD}}}
\def\om{\omega}
\def\sg{\sigma}
\def\g{\gamma}

\def\Sing{{\rm \text{Sing}}}
\def\Crit{{\rm \text{Crit}}}
\def\PS{{\rm \text{PS}}}
\def\P{\text{{\rm P}}}
\def\Comp{\text{{\rm Comp}}}
\def\H{{\rm H}}

\def\es{\emptyset}
\def\sbt{\subset}
\def\spt{\supset}
\def\sms{\setminus}
\def\bu{\bigcup}
\def\bi{\bigcap}
\def\({\big(}
\def\){\big)}

\newcommand{\h}{{\rm \text{h}}}
\newcommand{\ka}{\kappa}
\newcommand{\bpf}{\begin{proof}}
\newcommand{\epf}{\end{proof}}

\newcommand{\fr}{\noindent}
\newcommand{\spa}{\smallskip}

\newcommand{\Om}{\Omega}

\theoremstyle{definition}
\newtheorem{dfn}{Definition}[section]{\bf}{\rm}
\numberwithin{dfn}{section}
\newtheorem{rem}[dfn]{Remark}{\bf}{\rm}
{\bf}{\rm}

\theoremstyle{plain}
\newtheorem{lem}[dfn]{Lemma}{\bf}{\it}
\newtheorem{thm}[dfn]{Theorem}{\bf}{\it}
\newtheorem{prop}[dfn]{Proposition}{\bf}{\it}
{\bf}{\it}
\newtheorem{cor}[dfn]{Corollary}{\bf}{\it}
\newtheorem{fact}{Fact}{\bf}{\it}
\newtheorem{obs}[dfn]{Observation}{\bf}{\it}
\newtheorem{thma}{Theorem}{\bf}{\it}

\numberwithin{equation}{section}

\DeclareMathOperator{\diam}{diam}
\DeclareMathOperator{\Int}{Int}
\newcommand{\1}{{\mathbbm 1}}
\newcommand{\R}{\mathbb{R}}
\newcommand{\N}{\mathbb{N}}
\newcommand{\C}{{\mathbb C}}
\newcommand{\oc}{\widehat{\C}}
\def\sp{\medskip}
\def\b{\beta}
\makeatletter
\@namedef{subjclassname@2010}{%
  \textup{2010} Mathematics Subject Classification}
\makeatother

\begin{document}
\title[Thin annuli property, exponential distribution of return times]{Thin annuli property and exponential distribution of return times for~Weakly~Markov systems}
\thanks{The research of Mariusz Urba\'nski was supported in part by the NSF Grant DMS 1361677. Anna Zdunik was supported in part by the NCN grant 2014/13/B/ST1/04551}
\author{{\L}ukasz Pawelec}
\address{{\L}ukasz Pawelec, Department of Mathematics and Mathematical Economics, SGH Warsaw School of Economics, 
al.~Niepodleg\l{}o\'{s}ci~162, 02-554 Warszawa, Poland}
\email{LPawel@sgh.waw.pl}

\author{Mariusz Urba\'{n}ski} \address{Mariusz Urba\'{n}ski, University of North Texas, Department of Mathematics, 1155 Union Circle \#311430, Denton, TX 76203-5017,
USA}\email{Urbanski@unt.edu}

\author{Anna Zdunik}
\address{Anna Zdunik, Institute of Mathematics, University of Warsaw,
ul.~Banacha~2, 02-097 Warszawa, Poland}
\email{A.Zdunik@mimuw.edu.pl}
\begin{abstract}
We deal with the problem of asymptotic distribution of first return times to shrinking balls under iteration generated by a large general class of dynamical systems called weakly Markov. Our ultimate main result is that these distributions converge to the exponential law when the balls shrink to points. We apply this result to many classes of smooth dynamical systems that include conformal iterated function systems, rational functions on the Riemann sphere $\oc$, and transcendental meromorphic functions on the complex plane $\C$. We also apply them to expanding repellers and holomorphic endomorphisms of complex projective spaces. 


One of the key ingredients in our approach is to solve the  well known, in this field of mathematics, problem of appropriately estimating the measures of, suitably defined, large class of geometric annuli. We successfully do it. This problem is, in the existing literature, differently referred to by different authors; we call it the Thick Thin Annuli Property. 

Having this property established, we prove that for non--conformal systems the aforementioned distributions converge to the exponential one along  sets of radii whose relative Lebesgue measure converges fast to one.  

But this is not all. In the context of conformal iterated function systems, we establish the Full Thin Annuli Property, which gives the same estimates for all radii.  ln this way,  we solve   a long standing  problem.
As a result, we prove that the convergence to the exponential law holds along all radii for essentially all conformal iterated function systems and, with the help of the techniques of first return maps, for all aforementioned conformal dynamical systems.  
\end{abstract}
\keywords{return time, decay of correlations, exponential law}
\subjclass[2010]{ Primary 37B20, Secondary 37A25}

\maketitle

\tableofcontents
\section{Introduction}
In this paper we deal with asymptotic statistics of return times to shrinking objects that are formed by ordinary open balls with radii converging to zero. Let $(T,X,\mu,\rho)$ be a metric measure preserving dynamical system. By this we mean that $(X,\rho)$ is a metric space and $T\colon X\to X$ is a Borel measurable map preserving a Borel probability measure $\mu$ on $X$. 
Given a set $U\subset X$ and $x\in X$ define
\[
\tau_U(x):=\min\{n\ge 1:T^n(x)\in U\},
\]
and call it the (first) entry time to $U$.  When $x\in U$, this is called the first return time.

The modern study of return times was initiated in the early nineties by the seminal papers of M.~Boshernitzan \cite{Boshernitzan}, and D. S. Ornstein and B. Weiss \cite{OW}. They looked at return times to shrinking balls (Boshernitzan) or to decreasing cylinders in a symbol space (Ornstein, Weiss). These papers triggered a growing interest in the statistics of return times reflected in numerous publications on the subject. 

This lead, amongst others, to the notions of the \emph{recurrence dimensions} in \cite{BS} or \cite{STV}, the study of possible limiting distributions in \cite{KL}, speed of convergence, e.g. in \cite{todd}. One of the main avenues of study is the continuing search for systems for which the return (and entry) times distribution tends to the exponential law. By presenting the first general technique of proving such convergence, it was shown in \cite{HSV} that this should be the \emph{natural} limit. 

The early results focused on return times to cylinder sets. Passing to balls introduces some, qualitatively new, geometric flavour. It however ushers a considerable obstacle demanding to have a subtle upper estimates of measures of shrinking annuli. Virtually every proof of convergence to the exponential law requires that the measure of a thin annuli, of inner radii $r$ and outer radii $r+r^\kappa$, $\kappa>1$, are small compared to the measures of the balls of radii $r$. In all, known to us, papers such estimates are simply introduced as hypotheses, see for example the assumption (A4) from \cite{HW}, or the assumption (IV') from \cite{SRV}. This notion is also essential while proving the Poisson law as we may see in \cite{SP}, where the entire Appendix A in is devoted to comments on this problem.


This property is easily checked for measures equivalent to the Lebesgue measure or if the measure of every ball is bounded above by its radius raised to the power larger than $d-1$, where $d$ is the dimension of the ambient space.
It should be underlined that these measures were, essentially, the only ones for which the required aforementioned property was known. Consequently, the exponential limiting law for the return times all radii was known only for such measures. 

In the present paper we provide a fairly complete solution to the problem of estimating the measures of thin annuli, getting two results: firstly, see Theorem~\ref{t:thin-Intro}  and Section~\ref{sec:IFS}, we prove that it holds for conformal systems and all radii, and, secondly, see Theorem~\ref{Main_Thm_3}, for non-conformal ones, we show that it holds along a very large set of radii. Up to our knowledge, the present paper is the first one to tackle successfully the issue of proving upper estimates of measures of thin annuli.

\,\, In what follows we denote the open ball of radius $r>0$ centred at a point $x\in X$ either by both $B(x,r)$ or $B_r(x)$ depending on the context. Recall that given a measurable set $A$ with positive measure, we denote by $\mu_A$ the  measure on $A$, given by the formula
$$
\mu_A(F):=\frac{\mu(F)}{\mu(A)}
$$
where $F$ ranges over  measurable subsets of $A$.

Our main motivation and the main goal in this article is to identify a large  class of systems and large classes of families of positive radii $\cR=\{R_x\sbt (0,1]\}$, such that $0\in \overline R_x$, which are defined for $\mu$--a.e. $x\in X$, and for which the following properties hold:
\begin{equation}\label{eq:wykl_LM1}
\lim_{R_x\ni r\to 0}	\sup_{t\geq0}\left|\mu\left(\left\{z\in X\colon \tau_{B_r(x)}(z)>\frac{t}{\mu(B_r(x))}\right\}\right)-e^{-t}\right|= 0
\end{equation}
for $\mu$--a.e. $x\in X$, i.e. the distributions of the normalized first entry time converge to the exponential one law, and 
\begin{equation}\label{eq:wykl_LM2}
\lim_{R_x\ni r\to 0}	\sup_{t\geq0}\left|\mu_{B_r(x)}\left(\left\{z\in B_r(x)\colon \tau_{B_r(x)}(z)>\frac{t}{\mu(B_r(x))}\right\}\right)-e^{-t}\right|= 0
\end{equation}
for $\mu$--a.e. $x\in X$, i.e.
the distributions of the normalized first return time converge to the exponential one law. Formulas \eqref{eq:wykl_LM1} and \eqref{eq:wykl_LM2} are equivalent to saying that for every Borel set $F\sbt [0,+\infty)$ with boundary of Lebesgue measure zero, we have that, for $\mu$--a.e. $x\in X$, both of the following hold:
\begin{align}
\label{eq:wykl_LM1B}
\lim_{R_x\ni r\to 0}\mu\left(\left\{z\in X\colon \tau_{B_r(x)}(z){\mu(B_r(x))}\in F\right\}\right)&=\int_Fe^{-t}\,dt\\
\label{eq:wykl_LM2B}
\lim_{R_x\ni r\to 0}	\mu_{B_r(x)}\left(\left\{z\in B_r(x)\colon \tau_{B_r(x)}(z){\mu(B_r(x))}\in F\right\}\right)&=\int_Fe^{-t}\,dt.
\end{align}

\spa Our large class of measure preserving dynamical systems is that of \emph{Weakly Markov} ones defined, somewhat lengthily but naturally, in Section~\ref{sec:exp}. It is motivated by the class of loosely Markov systems introduced and explored in \cite{Urb2}. This class captures systems, not necessarily conformal, such as expanding repellers and holomorphic endomorphisms of complex projective spaces, but also conformal ones such as conformal graph directed Markov systems, conformal expanding repellers, rational functions of the Riemann sphere, and transcendental meromorphic functions. All this is described in detail in Section~\ref{sec:ex} devoted to examples. 

Having conformality in the system is not just to work in a more comfortable setting, but it does have seminal qualitative impact on the range of radii for which our main theorems hold. They do hold for all radii. 
Up to our best knowledge, this is the first time that for such general classes of systems and invariant measures the convergence to the exponential law is proved to hold along all radii.
\subsection{Subsets of radii}
In what follows Leb denotes Lebesgue measure. In this subsection $X$ is an arbitrary set. We now describe several natural classes of radii for which we will prove the aforementioned convergence to the exponential law. 
\begin{itemize}
\item The first class, called \emph{full}, denoted by $\cF$, contains all families $\{T_x:x\in X\}$, for which $T_x=(0,\eta_x]$ for some $\eta_x>0$.

\,

\item The next class, denoted by $\mathcal{AF}$, called \emph{almost full}, consists of all families $\{T_x:x\in X\}$ for which $\Leb(T_x\cap (0,\eta_x])=\eta_x$ for  some $\eta_x>0$.

\, 

\item The third class, denoted by $\mathcal{D}$ and called \emph{dense}, contains the families $\{T_x:x\in X\}$ satisfying 
\[
\lim_{r\to 0}\frac{\Leb(T_x\cap (0,r])}{r}=1 \mbox{\quad for all $x$}
\]
i.e., $0$ is the density point of $T_x$.

\, 

\item The fourth class, \emph{super dense} set denoted by $\mathcal{SD}$. It is composed of all families $\{T_x:x\in X\}$ for which,   for every $\alpha >0$
\[
\lim_{r\to 0}\frac{\left|\frac{\Leb(T_x\cap (0,r])}{r}-1\right|}{r^\alpha}=0 \mbox{\quad for all $x$.}
\]
\item Finally, the class class, called \emph{$\beta$--thick} (for $\beta>0$), and denoted by $\mathcal{\beta T}$ contains all families $\{T_x\}$ satisfying
\[
\lim_{r\to 0}\frac{\left|\frac{\Leb(T_x\cap (0,r])}{r}-1\right|}{r^{\ln^\beta(1/r)}}=0 \mbox{\quad for all $x$.}
\]
\end{itemize}
Trivially, for any $\beta>0$,
\[
\mathcal{F}\subset\mathcal{AF}\subset\mathcal{\beta F}\subset\mathcal{SD}\subset\mathcal{D}.
\]
\subsection{Main results}
\begin{thma}\label{wykl_LM2}
If $X$ is a Borel subset of $\R^d$ and $(T,X,\mu,\rho)$ is a Weakly Markov system, then for every $\beta>0$ there exists $\mathcal{\beta T}=\{T_x:x\in X\}$, a $\beta$--thick class of radii, such that for $\mu$--a.e. $x\in X$
\begin{equation}\label{eq:wykl_LM1C_1}
\lim_{T_x\ni r\to 0}	\sup_{t\geq0}\left|\mu\left(\Big\{z\in X\colon \tau_{B_r(x)}(z)>\frac{t}{\mu(B_r(x))}\Big\}\right)-e^{-t}\right|= 0,
\end{equation}
i.e. the distributions of the normalized first entry time converge to the exponential one law, and 
\begin{equation}\label{eq:wykl_LM2C_1}
\lim_{T_x\ni r\to 0}	\sup_{t\geq0}\left|\mu_{B_r(x)}\left(\Big\{z\in B_r(x)\colon \tau_{B_r(x)}(z)>\frac{t}{\mu(B_r(x))}\Big\}\right)-e^{-t}\right|= 0,
\end{equation}
i.e. the same convergence holds for the normalized first return time.
\end{thma}
\begin{rem}
In fact, as we prove, there are even larger classes of radii for which Theorem~\ref{wykl_LM2} holds. See Theorem~\ref{Main_Thm_2} and Remark~\ref{t1_2016_08_08}.
\end{rem}
We now introduce a crucial property of a measure, which we call the \emph{Thin Annuli Property}. We define and discuss it now in two steps. 

\begin{dfn}\label{def:subpoly}
A function $\kappa\colon (0,1]\to \R_+$ will be called \emph{subpolynomial} if it is monotone nonincreasing, and for every $\varepsilon>0$
\begin{equation}
\lim_{r\to 0} \kappa(r)r^\varepsilon = 0.
\end{equation}
\end{dfn}

\begin{rem}\label{r6_2016_09_06}
Subpolynomial functions include all positive constant functions and functions of the form $\kappa(r)=\alpha\ln^\beta(1/r)$, for $\alpha, \beta>0$.
\end{rem}

\begin{dfn}\label{d1_2016_06_16}
Let $\mu$ be a finite Borel measure on a metric space $X$. Let $\cR=\{R_x:x\in X\}$, be a class of radii defined $\mu$--a.e. in $X$. The measure $\mu$ is said to have the \emph{Thin Annuli Property} relative to $\cR$ if
for $\mu$--almost every $x\in X$ there exists a subpolynomial function $\kappa_x: (0,1]\to \R_+$ such that 
\begin{align}\label{eq:thanndef}
\lim_{R_x\ni r\to 0} \frac{\mu\left(B(x,r+r^{\kappa_x(r)})\setminus B(x,r)\right)}{\mu\left(B(x,r)\right)}=0.
\end{align}
We say that measure $\mu$ satisfies the \emph{Thick Thin Annuli Property} if for every $\beta>0$ it has the Thin Annuli Property with respect to some $\beta$--\emph{thick} class of radii.
We analogously define the \emph{Full Thin Annuli Property}.
\end{dfn}
The two main ingredients of the proof of Theorem~\ref{wykl_LM2}, and important results on their own, are the following.

\begin{thma}\label{Main_Thm_2}
Let $(T,X,\mu,\rho)$ be a Weakly Markov system. If the measure $\mu$ has the Thin Annuli Property relative to some class of radii $\cR$
, then both (\ref{eq:wykl_LM1}) and (\ref{eq:wykl_LM2}) hold, i.e. the distributions of the normalized first entry time and first return time converge to the exponential one law. 
\end{thma}

\begin{thma}\label{Main_Thm_3}
Every finite Borel measure $\mu$ in a Euclidean space $\R^d$, satisfies the Thick Thin Annuli Property.
\end{thma}

\begin{rem}\label{t1_2016_08_08}
In fact, Theorem \ref{Main_Thm_3} could be strengthened: see Theorem~\ref{t1_2016_09_02} along with  Theorem~\ref{thm:BSG}, Definition~\ref{def:subpoly} (and Remark~\ref{r6_2016_09_06}), and Remark~\ref{r5_2016_09_06}.
\end{rem}
\fr Of course, Theorem~\ref{wykl_LM2} is an immediate consequence of Theorems~\ref{Main_Thm_2} and~\ref{Main_Thm_3}.

The further natural question to ask is about the convergence to the exponential law along a full class of radii. Because of Theorem~\ref{Main_Thm_2} the answer would be positive if we had a Weakly Markov system whose measure has the Full Thin Annuli Property.

We have discovered that this property is satisfied for a large class of systems. The only additional requirement is for the system to be generated by a countable (either finite or infinite) alphabet conformal iterated function system (IFS). This gives rise, via suitable inducing schemes, to several classes of applications as it is shown in the last section. 

All the definitions appearing in the following theorems are in Subsection \ref{sec:IFS}. Here is the fourth main result, an achievement on its own.
\begin{thma}[see Theorem~\ref{t:thin}]\label{t:thin-Intro}
If $\cS =\{\phi_e\colon X\to X\}_{e\in E}$ is a conformal geometrically irreducible IFS, then for every $\mu\in\mathcal M_E$, a large class of measures containing many Gibbs/equilibrium measures of H\"{o}lder continuous summable potentials on the symbol space $E^\N$, the projection measure $\mu\circ\pi^{-1}$ on $J_\cS$ has the Full Thin Annuli Property. In fact, the following holds 
\begin{align*}
\lim_{r\to 0} \frac{\mu\circ\pi^{-1}\big(B(x,r+r^3)\sms B(x,r)\big)}{\mu\circ\pi^{-1}\left(B(x,r)\right)}=0\mbox{\qquad for \ $\mu\circ\pi^{-1}$--a.e. }	 x\in J_\cS.
\end{align*}
\end{thma}

\fr We should immediately emphasize that in this theorem the conformal IFS $\cS$ is not required to satisfy any kind of separation condition, nor even its weakest form known as the Open Set Condition. In other words, all kinds of overlaps are allowed. Also, the measures $\mu\in \mathcal M_E$ need not be Gibbs/equilibrium states nor even shift-invariant. These measures are just to satisfy two natural conditions formulated in Subsection~\ref{sec:IFS}. Theorem~\ref{t:thin-Intro} via Theorem~\ref{Main_Thm_2} leads to the convergence to the exponential distribution along all radii (full class) for all Weakly Markov systems generated by conformal IFS and Gibbs/equilibrium measures (now we do need them).

\begin{thma}[see Theorem~\ref{t1alr4}]\label{t1alr4-Intro}
Suppose that $\cS$ is a finitely irreducible and geometrically irreducible conformal IFS satisfying the Strong Open Set Condition. If  $f\colon E^\N\to\R$ is a summable H\"older continuous potential such that for some $\beta>0$
\begin{equation}\label{1_2016_02_22_1}
\sum_{e\in E}\exp\(\inf\(f|_{[e]}\)\)\|\phi_e'\|^{-\beta}<+\infty,
\end{equation}
then the measure--preserving dynamical system $\(T_\cS\colon\mathring{J}_\cS\to\mathring{J}_\cS,\hat\mu_f\)$ is Weakly Markov and satisfies the Full Thin Annuli Property. In consequence, the exponential one laws hold along all radii. 
\end{thma}
We highlight again that the above theorem is very general and allows us to prove, using a suitable inducing procedure, the exponential law for several naturally occurring classes of conformal systems, as seen in Section~\ref{sec:ex} devoted to applications and examples. 
 

\spa We end this introduction with a comment on the Weakly Markov systems. This concept captures and extends that of Loosely Markov systems of \cite{Urb2}, and of the earliest works on the subject such as \cite{STV}. One of the advantages of Weakly Markov systems is that no transfer operator is involved, and merely the exponential decay of correlations is assumed,  along with two other standard hypotheses.

\section{Convergence to Exponential Distribution for Weakly Markov Systems}\label{sec:exp}

In this section we do two things. First, we define the class of Weakly Markov systems and then we prove Theorem~\ref{Main_Thm_2}. We begin by
recalling the following standard definition:
\begin{dfn}\label{def:pointdim}
For a finite Borel measure $\mu$ on a metric space $X$, define the \emph{lower} and \emph{upper  pointwise dimensions}, denoted respectively by  $\underline{d}_{\mu}$ and $\overline{d}_{\mu}$, of the measure $\mu$ by 
\[\underline{d}_{\mu}(z) = \liminf_{r\to0} \frac{\ln\big(\mu(B_r(z))\big)}{\ln r},\hspace{3em} \overline{d}_{\mu}(z) = \limsup_{r\to0} \frac{\ln\big(\mu(B_r(z))\big)}{\ln r}.\] 
\end{dfn}

Passing to the next concept we need, given $\xi\in(0,1]$ denote by $\mathcal{H}^\xi(X)$ the vector space of all real--valued H\"older continuous functions on a metric space $(X,\rho)$ with exponent $\xi$, i.e. $f\in\mathcal{H}^\xi(X)$ if $f\colon X\to\R$ is bounded, continuous, and $v_\xi(f)<\infty$, where
\begin{equation}\nonumber
	v_\xi(f) := \inf\left\{H\ge 0 : \forall_{x,y\in X}\, |f(x)-f(y)|\leq H\rho^\xi(x,y))\right\}.
\end{equation}
The space $\mathcal{H}^\xi(X)$ is commonly endowed with the norm:
\begin{equation}\nonumber 
	||f||_\xi := ||f||_\infty + v_\xi(f),
\end{equation}
and then it becomes a Banach space.

Define further the first return of a set $U$ to itself under the map $T$ by 
\[
\tau(U):=\min_{x\in U}\tau_U(x).
\] 
\begin{dfn}\label{dfn:Weakly_Markov}
We will call a metric measure preserving dynamical system $(T,X,\mu,\rho)$, defined in Introduction, \emph{Weakly Markov}, if it satisfies the following conditions (i) to (iii): 
\begin{enumerate}[(i)]	
	\item \emph{Exponential Decay of Correlations}: There exists $\gamma \in (0,1)$ and $C>0$ (in general depending on $\xi$) such that for all $g\in \mathcal{H}^\xi$, all $f \in L^{\infty}_\mu$ and every $n\in \N$, we have 
\begin{equation}\label{eq:decaydef}
\left|\mu\left(f\circ T^n\cdot g\right)-\mu(g)\cdot\mu(f)\right| \leq C\gamma^n||g||_\xi\mu(|f|).
\end{equation}
	\item For $\mu$--a.e. $x\in X$, we have that
	$
0<\underline{d}_{\mu}(x)\leq \overline{d}_{\mu}(x)<+\infty.
$
\item \emph{No small returns}: $\displaystyle \liminf_{r\to 0} \frac{\tau\big(B_r(x)\big)}{-\ln(r)} > 0$ \;\; for $\mu$--a.e. $x\in X$.
\end{enumerate}
In addition, if measure $\mu$ also has the \emph{thin annuli property} relative to a family $\cR$ of radii, then we will call the system \emph{Weakly Markov with thin annuli} relative to $\cR$. If $\cR$ is thick (resp. full) we will call the system \emph{Weakly Markov with thick} (resp. \emph{full}) \emph{thin annuli}.
\end{dfn}\label{def:weakM}

\begin{rem}The \emph{no small returns} property has been proved to hold for many dynamical systems; e.g. those considered in \cite{STV}, and also, as it is easy to check, it holds for open transitive distance expanding maps and measures $\mu$ being Gibbs/equilibrium states of H\"older continuous potentials. For further information on this notion  see for example \cite{S} and the references therein.
\end{rem}

\begin{rem}For some specific systems the exponential distribution of the limit of return times has already been proved assuming only \emph{polynomial} decay of correlations. We however work in a very general setting, and our method, suitable for such generality, does need faster, in fact exponential, decay. 
\end{rem}

\begin{rem}\label{rem:part} 
As mentioned in the introduction, the second named author introduced the concept of \emph{Loosely Markov} systems in \cite{Urb2}. These systems are required to satisfy (ii), a stronger version of (i), and a Weak Partition Existence Condition, which implies (iii), as it was observed in \cite{Urb2}. Since we will also make use of this condition, we now formulate it below.
\end{rem}

\begin{dfn}\label{def:weak_partition} A system is said to satisfy the Weak Partition Existence Condition if there exists a countable partition $\alpha$ with entropy $\h_\mu(f,\alpha)>0$ and such that for $\mu$--a.e. $x\in X$ there exists $\chi(x)>0$ such that 
\begin{equation}\label{borel_cantelli}
B\(x,\exp(-\chi(x)n)\)\subset \alpha^n(x)
\end{equation}
for all integers $n\ge 0$ sufficiently large, where $\alpha^n:=\bigvee_{j=0}^{n-1}T^{-j}(\alpha)$ is the $n$-th refinement of 
 $\alpha$ under the action of $T$ and $\alpha^n(x)$ denotes the element of the partition containing $x$.
\end{dfn}
In order to prove Theorem~\ref{Main_Thm_2} we will apply two theorems from \cite{HSV}. 





\begin{proof}[Proof of Theorem \ref{Main_Thm_2}]
Recall that $(T,X,\rho,\mu)$ is a Weakly Markov system. Let us start with some notation; we follow \cite{HSV}. For a fixed set $U\sbt X$ let us define
\begin{align*}
c(k,U)&:=\mu_U\left(\tau>k\right)-\mu\left(\tau>k\right),\\
c(U)&:=\sup_{k\in\N}\left|c(k,U)\right|.
\end{align*}
The first result from \cite{HSV}, valid in a fairly abstract context, is this:

\begin{thm}\label{hsv1_LM} For a transformation $T\colon X\to X$, preserving a probability measure $\mu$ on $X$, the distributions of both the first return time and first entry time differ from the exponential law by an expression which converges to 0 if both $\mu(U)$ and $c(U)$ go to 0. More precisely, for entry time
\begin{equation}
	\sup_{t\geq0}\left|\mu\left(\left\{z\in X\colon \tau_{U}(z)>\frac{t}{\mu(U)}\right\}\right)-e^{-t}\right|\leq d(U),
\end{equation}
and also for return time
\begin{equation}
	\sup_{t\geq0}\left|\mu_U\left(\left\{z\in U \colon \tau_{U}(z)>\frac{t}{\mu(U)}\right\}\right)-e^{-t}\right|\leq d(U),
\end{equation}
where $d(U)=4\mu(U)+c(U)\big(1-\ln c(U)\big)$.
\end{thm}
\noindent The second theorem (also from \cite{HSV}) gives an estimate on the value of $c(U)$.
\begin{thm}\label{hsv2_LM}
With the transformation as above:
\[c(U) \leq \inf_{N\in\N} \left\{a_N(U) + b_N(U) + N \mu(U)\right\}, \mbox{\quad where}\]
\begin{align*}
a_N(U) &= \mu_U\left(\left\{\tau_U\leq N\right\}\right),\\
	b_N(U) &= \sup_{V\in \mathcal{B}}\left|\mu_U\big(T^{-N}V\big) - \mu\big(V\big)\right|=\sup_{V\in \mathcal{B}}\Big|\frac{\mu\big(U\cap T^{-N}V\big) -\mu(U) \mu\big(V\big)}{\mu(U)}\Big|,\\
	&\mbox{where $\mathcal{B}$ is the $\sigma$-algebra of Borel sets on $X$.} 
\end{align*}
\end{thm}
\begin{rem} Note that for a fixed set $U$ the number $a_N(U)$ grows to $1$ as $N\to+\infty$, whereas $b_N(U)$ tends to 0 (provided that the system has some mixing properties). The tricky part is to find a number $N$ such that $b_N$ has already become small, but $a_N$ and $N\cdot\mu(U)$   have not  grown too big.
\end{rem}
The proof of Theorem~\ref{Main_Thm_2} is a consequence of those two theorems and the following lemma, which is our main technical result in this section. 
\begin{lem}\label{hsv3_LM}
If a system $(T,X,\mu,\mathcal{B},\rho)$ is Weakly Markov with Thin Annuli Property relative to a class $\cR=\{R_x:x\in X\}$ of radii, then for $\mu$--almost all $x\in X$ and all radii $r>0$ there are integers $n_r(x)\ge 1$ such that
\[
\lim_{r\to 0}a_{n_r(x)}(B_r(x))
=\lim_{R_x\ni r\to 0}b_{n_r(x)}(B_r(x))
=\lim_{r\to 0}n_r(x)\cdot\mu(B_r(x))
=0
\] 
for $\mu$--almost all $x\in X$.
\end{lem}
\begin{proof} 
We will write $B_r$ instead of $B_r(x)$, when dependence on $x$ is clear. Put 
\[
n_r=n_r(x):=\mu(B_r)^{-\theta}.
\]
Obviously, if $\theta<1$ we get $n_r\cdot\mu(B_r)\to 0$ instantly. So it remains to find $\theta$ such that both $a_{n_r}$ and $b_{n_r}$ will tend to 0.

First, rewrite the \emph{no small returns} assumption: there exist a Borel set $V\subset X$ of full $\mu$ measure and two measurable functions $\chi(x),\, \rho_1(x)$, both positive $\mu$--a.e., such that 
\begin{equation}\label{eq:nosmret}
B_r(x)\cap T^{-k}\big(B_r(x)\big)=\emptyset 
\end{equation} 
for all $x\in V$, all radii $0<r<\rho_1(x)$ and all integers $1\le k\leq \chi(x)\ln(1/r)$.

Secondly, the assumptions imposed on pointwise dimension imply that there exists a set $W\subset V\subset X$, of full measure $\mu$, such that for all $x\in W$ 
\begin{equation}\label{pd_LM}
	r^{2\overline{d}_{\mu}(x)}\leq \mu\left(B_r(x)\right)\leq r^{\underline{d}_{\mu}(x)/2},
\end{equation}
for all radii $0<r<\rho_2(x)$ with a certain measurable, positive $\mu$--a.e. function $\rho_2\le \rho_1$.
Now let us define a family of Lipschitz continuous functions approximating a characteristic function on a ball; depending on:  radius $r >0$, real number $\alpha>0$, and $x\in X$, which will vary in the sequel; particularly, we will utilize various choices of $\alpha>0$. We define the auxiliary functions
\begin{equation}\nonumber
\phi_{r}^{(\alpha)}(t)
:= \left\{\begin{array}{ll} 1 
& \mbox{for } 0\leq t \leq r \\ r^{-\alpha}(r+r^{\alpha}-t) 
& \mbox{for } r\leq t 
\leq r+r^{\alpha}\\ 0 
& \mbox{for } t\geq r+r^{\alpha}\end{array}\right..
\end{equation}
The approximating functions are 
\[
g_{r,x}^{(\alpha)}(z):=\phi_r^{(\alpha)}(\rho(z,x)).
\]
The Lipschitz constant of $g_{r,x}^{(\alpha)}$ is bounded from above by $\displaystyle r^{-\alpha}$ as metric $\rho$ is 1--Lipschitz. In particular their H\"older norms (needed in the definition of exponential decay of correlations) are bounded from above by (take $\xi=1$)
$$
||g_{r,x}^{(\alpha)}||_\xi \leq 1+r^{-\alpha}\approx r^{-\alpha}
$$ 
for all $r>0$ sufficiently small. Fix $x\in W\subset V$ and small $r>0$. Put $g_r:=g_{r,x}^{(\alpha)}$ and put $f_r:=\1_{B_r}$.
Note that $f_r\leq g_r$. Recall the definition  of $a_N$.
\begin{equation*}
	a_N(B_r)=\mu_{B_r}\left(\tau_{B_r}\leq N\right)=\mu_{B_r}\left(\bigcup_{n=1}^N T^{-n}(B_r)\right)\leq \sum_{n=1}^N \frac{\mu\left(B_r\cap T^{-n}(B_r)\right)}{\mu(B_r)}.
\end{equation*}
As $x\in V$ equation \eqref{eq:nosmret} gives that some first intersections are empty. Put $\chi:=\chi(x)$.
\begin{equation*}a_N(B_r)
\leq \sum_{n=-\chi\ln(r)}^N \frac{\mu\left(B_r\cap T^{-n}(B_r)\right)}{\mu(B_r)}.
\end{equation*}
The assumption (\ref{eq:decaydef}) on decay of correlations gives
\begin{align*}
	\mu\left(B_r\cap T^{-n}(B_r)\right)
	&=\mu\left(f_r\circ T^n\cdot f_r \right)
	\leq \mu\left(f_r\circ T^n\cdot g_r\right)\\
	&\leq \mu(g_r)\cdot\mu(f_r) + C\gamma^n||g_r||_\xi\mu(f_r)\\
	&\leq \mu(f_r)\left(\mu(g_r) + C\gamma^n r^{-\alpha}\right).
\end{align*}
This allows us to rewrite the estimate on $a_N$ and later to bound the sum's elements as simply as possible in the following way
\begin{equation}\label{1_2015_12_03}
\begin{aligned}
a_N(B_r)
&\leq \sum_{n=-\chi\ln(r)}^N \left(\mu(g_r) + C\gamma^n r^{-\alpha}\right)	
 \leq N\mu(g_r)+ Cr^{-\alpha}\sum_{n=-\chi\ln(r)}^{+\infty} \gamma^n \\
&= N\mu(g_r)+ \frac{C}{1-\gamma}r^{-\alpha}\gamma^{-\chi\ln(r)}
=N\mu(g_r)+ Dr^{-\alpha-\chi\ln(\gamma)}.
\end{aligned}
\end{equation}
Now specify $\alpha>0$ to be in $(0,1]$. Using (\ref{pd_LM}) we estimate as follows.
\begin{equation}\label{2_2015_12_03}
\mu\(g_r\)
\leq \mu\left(B(x,r+r^\alpha)\right)
\leq \mu\left(B(x,2r^\alpha)\right) 
\leq 2^{\underline{d}_{\mu}(x)/2}r^{\alpha\underline{d}_{\mu}(x)/2}.
\end{equation}
Take $\theta>0$ as small as needed in the course of the proof and fix 
 $N=n_r=\mu(B_r)^{-\theta}$.  Inserting \eqref{2_2015_12_03} into \eqref{1_2015_12_03}, and using \eqref{pd_LM} again, we get
\begin{equation}\label{1_2016_08_09}
\begin{aligned}
a_{n_r}(B_r)
&\leq \mu(B_r)^{-\theta}\cdot 2^{\underline{d}_{\mu}(x)/2}r^{\alpha\underline{d}_{\mu}(x)/2}+ Dr^{-\alpha-\chi\ln(\gamma)}\\
&\leq E r^{-2\theta \overline{d}_{\mu}(x)}r^{\alpha\underline{d}_{\mu}(x)/2}+ Dr^{-\alpha-\chi\ln(\gamma)},
\end{aligned}
\end{equation}
with some positive constants $D$ and $E$. Restrict further the choice of $\alpha>0$ so that $\alpha<-\chi\ln(\gamma)$. Then fix any $\theta>0$ so small that $2\theta \overline{d}_{\mu}(x)<\alpha\underline{d}_{\mu}(x)/2$. With these specifications both exponents of powers of $r$ in formula \eqref{1_2016_08_09} are positive; whence we arrive at the conclusion that
\begin{equation}
\lim_{r\to 0} a_{n_r}(B_r) = 0.
\end{equation}
Note that our reasoning leading to this formula did not require any kind of the thin annuli property at all.

\spa Now we turn to the task of estimating $b_{n_r}(B_r(x))$. For this we do need and we do use the Thin Annuli Property relative to $\cR$. 

Let $\kappa_x:(0,1]\to (0,+\infty)$ 
be the subpolynomial function resulting from the Thin Annuli Property of the system $(T, X,\mathcal{B},\rho)$ relative to $\cR$. The point $x\in W$ is as above and also respecting formula (\ref{eq:thanndef}) of Definition~\ref{d1_2016_06_16}. Put
$g_r:= g_{r,x}^{(\kappa_x(r))}$ and fix a Borel set $H$. Then
\begin{align*}
\big|\mu\big(B_r&\cap  T^{-n_r}(H)\big) -\mu(B_r) \mu\big(H\big)\big|=\\
&= \big|\mu\big(\1_H\circ T^{n_r}\cdot f_r\big) -\mu(\1_H) \mu(f_r)\big|\le \\
&\le \left|\mu\big(\1_H\circ T^{n_r}\cdot f_r\big) -\mu\big(\1_H\circ T^{n_r}\cdot g_r\big)\right|+ \\
& \  \  \ +\left|\mu\big(\1_H\circ T^{n_r}\cdot g_r\big) -\mu(\1_H) \mu(g_r)\right|
  +\big|\mu(\1_H) \mu(g_r) -\mu(\1_H) \mu(f_r)\big|.
\end{align*}
So $\mu(B_r)b_{n_r}(B_r)$ is bounded by the supremum (over all Borel sets $H\subset X$) of the sum of the three terms above.  

The third expression bounding $b_{n_r}(B_r)$ is estimated easily:
\begin{align*}
	\mu(B_r)^{-1}\big|\mu(\1_H) \mu(g_r) &-\mu(\1_H) \mu(f_r)\big|\le \\
&\leq \mu(B_r)^{-1}\big( \mu(g_r) - \mu(f_r)\big)\ \\
	&\leq \mu(B_r)^{-1} \big( \mu(B(x,r+r^{\kappa_x(r)}) - \mu(B(x,r)\big)\\
	&= \frac{\mu\left(B(x,r+r^{\kappa_x(r)})\setminus B(x,r)\right)}{\mu(B_r)}.
\end{align*}
This tends to $0$ as $R_x\ni r\to 0$ because of the Thin Annuli Property relative to $\cR$ assumed to hold.
The first term is bounded in the same way since 
\begin{align*}
	\left|\mu\big(\1_H\circ T^{n_r}\cdot f_r\big) -\mu\big(\1_H\circ T^{n_r}\cdot g_r\big)\right| \leq  \mu(g_r) - \mu(f_r).
\end{align*}
Dealing with the second term we use the exponential decay of correlations:
\begin{align*}
\left|\mu\big(\1_H\circ T^{n_r}\cdot g_r\big) -\mu(\1_H) \mu(g_r)\right|
\leq C\gamma^n_r r^{-\kappa_x(r)}\mu(\1_H)
\leq C\gamma^n_r r^{-\kappa_x(r)}.
\end{align*}
The pointwise dimensions formula (\ref{pd_LM}) gives $n_r = \mu(B_r)^{-\theta} \geq r^{-\theta\underline{d}_{\mu}(x)/2}$ and 
\begin{align*}
	\mu(B_r)^{-1}\Big|\mu\big(\1_H\circ T^{n_r}\cdot g_r\big) -\mu(\1_H) \mu&(g_r)\Big|
	\leq C r^{-\kappa_x(r) - 2\overline{d}_\mu(x)}\gamma^{r^{-\theta\underline{d}_{\mu}(x)/2}}\\
	&=Ce^{-\kappa_x(r)\ln(r)- 2\overline{d}_\mu(x)\ln(r) + r^{-\theta\underline{d}_{\mu}(x)/2} \ln(\gamma)}.
\end{align*}
The last term in this formula converges to zero as $r \to 0$ once we know that
\begin{equation}\label{eq:spap}
	\lim_{r\to 0} \kappa_x(r)\ln(r)r^{\theta\underline{d}_{\mu}(x)/2}= 0,
\end{equation}
which indeed holds because $\kappa_x$ is a subpolynomial function. Thus,  
\begin{equation}
	\lim_{\cR_x\ni r\to 0} b_{n_r}(B_r) = 0.
\end{equation}
and this ends the proof of Lemma~\ref{hsv3_LM}.
\end{proof}
\fr The proof of Theorem~\ref{Main_Thm_2} is complete.
\end{proof}
\section{The Thin Annuli Property}

\subsection{Thick Thin Annuli Property holds for  Finite Borel Measures in $\mathbb R^d$.}\label{sec:sub}  

\

\noindent Our main result in this subsection is Theorem~\ref{Main_Thm_3}.
 We will need several technical auxiliary results, one of which, Theorem~\ref{t1_2016_09_02} is of high generality, interesting in itself, and entails Theorem~\ref{Main_Thm_3}. The following well-known result comes from \cite{BS}.
\begin{prop}\label{prop2:BS}
Any Borel probability measure on $\R^d$ is \emph{weakly diametrically regular}, i.e. 
for $\mu$--almost every $x\in \R^d$ and every $\varepsilon>0$ there exists $\delta>0$ such that for all $0<r<\delta$
\begin{equation}
\mu(B(x,2r))\leq \mu(B(x,r))r^{-\varepsilon}.	
\end{equation}
\end{prop}
\noindent We prove the following remarkable strengthening of this proposition.

\begin{thm}\label{thm:BSs}
Assume that $\mu$ is a Borel probability measure on $\R^d$ and fix any $\varepsilon >0$. Then for $\mu$--a.e. $x\in\R^d$ and every sufficiently small $r>0$ (i.e. $0<r\leq \delta(x)$ and $\delta(x)>0$ $\mu$--a.e.) we have 
\begin{equation}\label{eq:thmBS_1}
	\mu\left(B(x,2r)\right)\leq \log_2^{2+\varepsilon}(1/r) \mu(B(x,r)).
\end{equation}
Moreover, if $s$ and $r$ are such that $\left[-\log_2(s)\right]=\left[-\log_2(r)\right]$ (i.e. for some $k$ we have $2^{-k-1} <r,s \leq 2^{-k}$), then
\begin{equation}
	\left[-\log_2(r)\right]^{-1-\varepsilon} \mu\left(B(x,r)\right) \leq \mu\left(B(x,s)\right)\leq \left[-\log_2(r)\right]^{1+\varepsilon} \mu\left(B(x,r)\right).
\end{equation}
\end{thm}

\begin{proof}
Fix $(\alpha_n)_{n=1}^\infty$, a sequence of positive numbers converging to zero and define \emph{bad} sets
\begin{equation}
	Z_n:= \{x: \mu(B(x,2^{-n}))\cdot \alpha_n >  \mu(B(x,2^{-n-1}))\}.
\end{equation}
We will show that for every $n\ge 1$ we have that
$$
\mu(Z_n)\leq M_d\alpha_n,
$$
where $M_d$ is the constant resulting from Besicovitch's Covering Theorem. Indeed, by virtue of this theorem we can cover $Z_n$ by balls $B(z_i,2^{-n-1})$, $i\in I$, centred at the set $Z_n$, in such a way that this covering has multiplicity bounded above by $M_d$. This leads to the estimate
\begin{equation}
	\mu(Z_n) \leq \sum_{i\in I} \mu(B(z_i,2^{-n-1})) < \sum_{i\in I}  \alpha_n \mu(B(z_i,2^{-n})) \leq \alpha_n M_d \mu (\R^d) \leq \alpha_n M_d.
\end{equation}
So, if in addition, $\sum_n \alpha_n <+\infty$, then by Borel-Cantelli Lemma, $\mu$--a.e. $x\in\R^d$ belongs only to finitely many sets $Z_n$.

Now consider
$$
\alpha_n:=n^{-1-\varepsilon/2}.
$$ 
Of course $\sum_n \alpha_n <+\infty$ and let $x\in\R^d$ belong to the full measure set given by the Borel-Cantelli Lemma. This means that there exists $N=N(x)$ such that for every $n\geq N$, we have
\begin{equation}
	\mu(B(x,2^{-n})) \leq n^{1+\varepsilon/2} \mu(B(x,2^{-n-1})).
\end{equation}
Take $0<r\leq 2^{-N-1}$ and put $k=\left[-\log_2(r)\right]$, i.e. $2^{-k-1} <r \leq 2^{-k}$. Then
\begin{align}\label{eq1}
\nonumber	\mu(B(x,2r))&\leq \mu(B(x,2^{-k+1}))\leq (k-1)^{1+\varepsilon/2}\mu(B(x,2^{-k}))\\ 
	&\leq k^{1+\varepsilon/2}(k-1)^{1+\varepsilon/2}\mu(B(x,2^{-k-1}))\leq k^{2+\varepsilon} \mu(B(x,r))  \\
\nonumber	&\leq \left[-\log_2(r)\right]^{2+\varepsilon}\mu(B(x,r)) \leq \log_2^{2+\varepsilon}\left(\frac{1}{r}\right)\mu(B(x,r)). \qedhere
\end{align}
\end{proof}
\begin{rem}\label{r1_2016_08_15}
By taking a different convergent series, e.g. $\frac{1}{n\log^2(n)}$ as $\alpha_n$ we could improve the above estimate (and therefore in Cor. \ref{c1_2016_09_06}) to e.g.
\[
	\mu\left(B(x,2r)\right)\leq \log^2_2(1/r) \log^{2+\varepsilon}(\log(1/r)) \mu(B(x,r)).
\]
\end{rem}
Motivated by Theorem~\ref{thm:BSs} and Remark~\ref{r1_2016_08_15} we introduce the following.
\begin{dfn}\label{d2_2016_15}
A non--increasing function $G\colon(0,+\infty)\to (1,+\infty)$ satisfying
\begin{equation}\label{1_2016_09_02}
G(r/2)\le \gamma G(r)
\end{equation}
with some $\gamma\in[1,2)$ and all $r>0$ small enough, is called a \emph{doubling bound} for a Borel probability measure $\mu$ on $\R^d$ if for $\mu$--a.e. $x\in\R^d$, all sufficiently small $r>0$ (i.e. $0<r\leq \delta(x)$ and $\delta(x)>0$ $\mu$--a.e.), we have that 
\begin{equation}\label{eq:thmBS_2}
\mu\left(B(x,2r)\right)\leq G(r)\mu(B(x,r)).
\end{equation}
\end{dfn}
\noindent With this definition Theorem~\ref{thm:BSs} can be reformulated as follows. 
\begin{thm}\label{thm:BSG}
For every $\varepsilon>0$ the function
\[
(0,+\infty)\ni r\longmapsto \max\big\{0,\log_2^{2+\varepsilon}(1/r)\big\},
\]
(in fact any function of Remark~\ref{r1_2016_08_15}) is a \emph{doubling bound} for any Borel probability measure $\mu$ on $\R^d$.
\end{thm}
Definition~\ref{d2_2016_15} and Theorem~\ref{thm:BSG} will lead us to the following crucial technical estimate on the measures of annuli. 
\begin{lem}\label{thmmiar2*}
Let $\kappa_x\colon (0,1]\to\R$ 
be subpolynomial functions such that
\begin{equation}\label{1_2016_08_12}
\underline\kappa_x:=\inf_{r\in(0,1]}\kappa_x(r)>1 \mbox{\qquad for every $x\in \R^d$.}
\end{equation}
If $\mu$ is a Borel probability measure on $X=\R^d$, then for $\mu$--a.e. $x\in X$ and every  $A>0$ the set of those radii $r>0$ for which
\begin{align}\label{eq2:thmmiar2}
\frac{\mu\left(B(x,r+r^{\kappa_x(r)})\setminus B(x,r)\right)}{\mu\left(B(x,r)\right)}>A	
\end{align}
has zero density at the point $r=0$.
In other words, if we denote by $Z_x(A)$ the set of all radii $r>0$ that satisfy \eqref{eq2:thmmiar2}, then
\begin{equation}\label{3_2016_09_02}
\lim_{r\to 0} \frac{l(Z_x(A)\cap [0,r])}{l([0,r])} = 0, \ \mbox{where $l$ is Lebesgue measure on $\R$.}
\end{equation}
Moreover, let $G$ be  a doubling bound almost everywhere for $\mu$. Then the following, more precise estimate holds:  
\begin{equation}\label{dokladne_szacowanie}
l(Z_x(A)\cap [0,r])  
\leq
\frac{2}{\left(1-\frac{\gamma}2\right)\ln(1+A)}r^{\kappa_x(r)}\ln G(r). 
\end{equation}
\end{lem}
\begin{proof}
The first observation is that \eqref{3_2016_09_02} follows from \eqref{dokladne_szacowanie}. Indeed, it suffices to take 
$G(r)=r^{-\alpha}$ for some $\alpha\in (0,1)$.
We are therefore to prove \eqref{dokladne_szacowanie} only. We do it now.

$G$ is a doubling bound, so there exists a Borel set $Y\sbt \R^d$ with $\mu(Y)=1$ such that for every $x\in Y$ there exists $\delta_x>0$ such that for all $r\in (0,2\delta_x)$
\begin{equation}\label{eq2:BS}
\mu\big(B(x,2r)\big)\leq G(r)\mu\big(B(x,r)\big).	
\end{equation}
Fix any $x\in Y$. Then fix $r\in(0,\delta_x)$ and $\eta>0$. There exist an integer $n\ge 1$, and a sequence of $n$ radii $r_j\in(0,\delta_x)\cap Z_x(A)$, $j=1,2\ldots, n$ such that 
\begin{equation}
	r\le r_1 < r_1+r_1^{\kappa_x(r_1)} < r_2 < r_2+r_2^{\kappa_x(r_2)} < r_3 < \cdots < r_n + r_n^{\kappa_x(r_n)} \leq 2r
\end{equation}
and
\begin{equation}\label{5_2016_08_15}
l\left(\Big(Z_x(A)\cap [r,2r]\Big)\sms \bu_{j=1}^n\Big[r_j,r_j+r_j^{\kappa_x(r_j)}\Big)\right)\le \eta.
\end{equation}
In particular the annuli defined by radii $r_j$ do not intersect. Since $r_j\in Z_x(A)$ for all $j=1,2\ldots, n$, for any $1\leq p \leq n$, we have that 
\begin{equation}
\frac{\mu\left(B\(x,r_p+r_p^{\kappa_x(r_p)}\)\right)}{\mu\left(B(x,r_p)\right)}>1+A.
\end{equation}
Using this estimate $n$ times we arrive at
\begin{equation}
\begin{aligned}
\mu\left(B(x,r)\right)
\le \mu\left(B(x,r_1)\right) 
&\leq \frac{\mu\left(B\(x,r_1+r_1^{\kappa_x(r_1)}\)\right)}{1+A} \leq  \frac{\mu\left(B(x,r_2)\right)}{1+A} \leq \cdots\\ \nonumber
\cdots &\leq \frac{\mu\left(B(x,r_n)\right)}{(1+A)^n} \leq 	\frac{\mu\left(B(x,2r)\right)}{(1+A)^n}.
\end{aligned}
\end{equation}
Applying further (\ref{eq2:BS}) yields
$$
\mu(B(x,r))\leq \frac{\mu(B(x,r))G(r)}{(1+A)^n}.
$$
This shows that
\begin{equation}
G(r) \geq (1+A)^n, \ \mbox{\ giving the estimate: \hskip0.3cm} n\leq \frac{\ln G(r)}{\ln(1+A)}.
\end{equation}
Now divide the interval $[r,2r)$ into subintervals of length $(2r)^{\kappa_x(2r)}$. Define: 
\[
I_1=\left[r, r+(2r)^{\kappa_x(2r)}\right), \;\ldots\;, I_k=\left[r+(k-1)(2r)^{\kappa_x(2r)}, r+k(2r)^{\kappa_x(2r)}\right)\ldots
\] 
for all $k\ge 1$ until $(k+1)(2r)^{\kappa_x(r)}\geq 2r$.

Observe that $r_p^{\kappa_x(r_p)}\le (2r)^{\kappa_x(2r)}$, since $r_p\le 2r$ and the function $\kappa_x$ is nonincreasing. So, if $r_p\in I_k$ then the interval $\big[r_p,r_p + r_p^{\kappa_x(r_p)}\)$ is contained in $I_k\cup I_{k+1}$.


Thus, the union 
\[
\bu_{j=1}^n\left[r_j,r_j+r_j^{\kappa_x(r_j)}\right)
\]
is contained in a union of at most $2n\le \frac{2\ln G(r)}{\ln(1+A)}$ intervals of the form $I_k$. So
\begin{equation}
l\bigg(\bu_{j=1}^n\big[r_j,r_j+r_j^{\kappa_x(r_j)})\bigg)	
\leq (2r)^{\kappa_x(2r)} \cdot \frac{2\ln G(r)}{\ln(1+A)} 
= \frac{2}{\ln(1+A)}(2r)^{\kappa_x(2r)}\ln G(r).
\end{equation}
Along with \eqref{5_2016_08_15}, this gives
$$
l\(Z_x(A)\cap [r,2r]\)\le \eta +\frac{2}{\ln(1+A)}(2r)^{\kappa_x(2r)}\ln G(r).
$$ 
Since $\eta>0$ was arbitrary, this in turn gives
$$
l\(Z_x(A)\cap [r,2r]\)\le \frac{2}{\ln(1+A)}(2r)^{\kappa_x(r)}\ln G(r).
$$
By summing this estimate and recalling that the function $\kappa_x$ is monotone decreasing while the function $G$ satisfies \eqref{1_2016_09_02}, we get
\begin{align*}
\nonumber l(Z_x(A)\cap [0,r]) 
&\leq \sum_{j=1}^\infty l\left(Z_x(A)\cap\left[\frac{r}{2^j},
        \frac{r}{2^{j-1}}\right]\right) \\&
        \leq \frac{2}{\ln(1+A)}\sum_{j=1}^\infty \left(\frac{r}{2^{j-1}}\right)^{\kappa_x(r/2^{j-1})}
     \ln G\left(\frac{r}{2^j}\right)\\
&\leq\frac{2}{\ln(1+A)}r^{\kappa_x(r)}\ln G(r)\sum_{j=1}^\infty(\gamma/2)^{j-1} \\&=\frac{2}{\left(1-\frac{\gamma}2\right)\ln(1+A)}r^{\kappa_x(r)}\ln G(r). \qedhere
\end{align*}	
\end{proof}
As a consequence of this lemma we get the following first main result of this section.
\begin{thm}\label{t1_2016_09_02}
Let $g\colon\R_+\to \R_+$ be a function such that $\lim_{r\to 0}g(r)=+\infty$
and for every $\alpha>0$, every $s>0$ sufficiently small, and every $0<r\le s$
$$
\frac{g(r)}{g(s)}\le \left(\frac{s}{r}\right)^\alpha.
$$
Let $\mu$ be a Borel probability measure on $X=\R^d$ and let $G$ be a doubling bound almost everywhere for $\mu$. For every $x\in \R^d$ let $\kappa_x:(0,1]\to (1,+\infty)$ be a subpolynomial function such that
\begin{equation}\label{5_2016_08_12}
\underline\kappa_x:=\inf_{r\in(0,1)}\kappa_x\(r)>1.
\end{equation}
Then the measure $\mu$ has the Thin Annuli Property with respect to some class of radii $\cR=\{R_x\}_{x\in X}$ for which
\begin{equation}
\lim_{R_x\ni r\to 0}\,	\frac{\left|\frac{l(R_x\cap(0,r])}{r}-1\right|}{g(r)r^{\kappa_x(r)-1}\ln G(r)}=0.
\end{equation}
In addition, the subpolynomial functions witnessing this Thin Annuli Property are just the functions $\kappa_x$ introduced above in the hypotheses. 
\end{thm}

\begin{proof}
We first shall prove the following.

\spa\fr{\bf Claim~$1^0$:} There exists a constant $Q\ge 1$ such that
$$
g(r)r^{\ka_x(r)}\ln G(r)\le Qg(s)s^{\ka_x(s)}\ln G(s)
$$
for every $s>0$ sufficiently small and every $0<r\le s$.
\begin{proof}[Proof of Claim $1^0$]
The formula of this claim is equivalent to the following:
$$
\frac{g(r)}{g(s)}\cdot \frac{\ln G(r)}{\ln G(s)}
\le Q\frac{s^{\ka_x(s)}}{r^{\ka_x(r)}}.
$$
Since the function $\ka_x$ is monotone decreasing, we have that
$$
\left(\frac{s}{r}\right)^{\underline\ka_x}
\le \left(\frac{s}{r}\right)^{\ka_x(s)}
\le \frac{s^{\ka_x(s)}}{r^{\ka_x(r)}}.
$$
It therefore suffices to show that
$$
\frac{g(r)}{g(s)}\cdot \frac{\ln G(r)}{\ln G(s)}
\le Q\left(\frac{s}{r}\right)^{\underline\ka_x}.
$$
And in order to have this it suffices to know that

\begin{equation}\label{1alr11}
\frac{g(r)}{g(s)}\le \left(\frac{s}{r}\right)^{\underline\ka_x/2} \mbox{\qquad and\qquad} \frac{\ln G(r)}{\ln G(s)}\le Q\left(\frac{s}{r}\right)^{\underline\ka_x/2}.
\end{equation}
The former follows directly (for $s>0$ small enough) from our hypotheses while for proving the latter fix a unique integer $k\ge 0$ such that $2^kr\le s$ and $2^{k+1}r>s$. Then
$$
G(r)\le \gamma^{k+1}G(2^{k+1}r)\le \gamma^{k+1} G(s).
$$
Therefore $\ln G(r)\le (k+1)\ln\gamma+\ln G(s)$. Hence
$$
\frac{\ln G(r)}{\ln G(s)}\le 1+\ln\gamma\frac{k+1}{\ln G(s)}.
$$
Thus in order to have \eqref{1alr11} it suffices to know that
$$
1+\ln\gamma\frac{k+1}{\ln G(s)}\le Q\cdot 2^{\frac12\underline\ka_xk}.
$$
But as $\inf\ln G>0$, this inequality clearly holds for a sufficiently large constant $Q\ge 1$, all integers $k\ge 0$ and all $s>0$ sufficiently small. The claim is proved.
\end{proof}
Passing to the actual proof of Theorem~\ref{t1_2016_09_02}, we note that by Lemma~\ref{thmmiar2*},  the estimate \eqref{dokladne_szacowanie}, and by the fact that $g(r)\to\infty$ as $r\to0$,  there exists $(r_n)_{n=1}^\infty$, a decreasing sequence of positive radii converging to $0$ such that
\begin{equation}\label{6_2016_09_02}
l(Z_x(1/n)\cap [0,r])\le 2^{-n}Q^{-1}g(r)r^{\kappa_x(r)}\ln G(r)
\end{equation}
for all integers $n\ge 1$ and all radii $r\in(0,r_n]$. For $x\in X$ define
\[
Z_x:=\bu_{n=1}^\infty Z_x(1/n)\cap (r_{n+1},r_n]
\]
and then $R_x:=(0,1)\sms Z_x$.
For every $r\in(0,r_1]$ let $n=n_r\ge 1$ be the unique integer such that $r_{n+1}<r\le r_n$. Using Claim~$1^0$, we then estimate
$$
\begin{aligned}
l(Z_x\cap (0,r]&)
=\sum_{k=n+1}^\infty l\(Z_x\cap (r_{k+1},r_k]\)+l\(Z_x\cap(r_{n+1},r]\) \\
&=\sum_{k=n+1}^\infty l\(Z_x(1/k)\cap (r_{k+1},r_k]\)+l\(Z_x(1/n)\cap(r_{n+1},r]\) \\
&\le \sum_{k=n+1}^\infty 2^{-k}Q^{-1}g(r_k)r_k^{\kappa_x(r_k)}\ln G(r_k)+
       2^{-n}Q^{-1}g(r)r^{\kappa_x(r)}\ln G(r) \\
&\le Q\sum_{k=n+1}^\infty 2^{-k}Q^{-1}g(r)r^{\kappa_x(r)}\ln G(r)+
      2^{-n}g(r)r^{\kappa_x(r)}\ln G(r) \\
&=\sum_{k=n_r}^\infty2^{-k}g(r)r^{\kappa_x(r)}\ln G(r) 
= 2^{-n_r+1}g(r)r^{\kappa_x(r)}\ln G(r).
\end{aligned}
$$
Therefore, since $\lim_{r\to 0}n_r=+\infty$, we get that
\[
\lim_{r\to 0}\frac{l(Z_x\cap (0,r])}{g(r)r^{\kappa_x(r)}\ln G(r)}
\le \lim_{r\to 0}2^{-n_r+1}
=0.\qedhere
\]
\end{proof}
\begin{rem}\label{r5_2016_09_06}
Note that any iterate of the logarithmic function,  
$$
g(r)=\ln^{(k)}(1/r), \quad  k\in \N,
$$
satisfies the hypotheses of Theorem~\ref{t1_2016_09_02}.
\end{rem}
Taking $g(r)=\ln\ln(\frac{1}{r})$  and, using also Theorem ~\ref{thm:BSG}, the functions $G(r)$ and $\kappa_x(r)$, being respectively of the form $r\mapsto\log_2^{2+\varepsilon}(1/r)$ and  $r\mapsto\ln^\beta(1/r)+1$, we get the following. 
\begin{cor}\label{c1_2016_09_06}
Every finite Borel measure $\mu$ on $\R^d$ for every $\beta>0$ has the  Thin Annuli Property with respect to some class of radii  $\cR(\beta)=\{R_x(\beta)\}_{x\in X}$ satisfying
\begin{equation}
\lim_{R_x(\beta)\ni r\to 0}\, \frac{\left|\frac{l(R_x(\beta)\cap(0,r])}{r}-1\right|}{r^{\ln^\beta(1/r)}\ln \ln(1/r)}=0.
\end{equation}
\end{cor}
\fr Ending this part we observe that this corollary directly entails Theorem~\ref{Main_Thm_3}.

\subsection{Conformal Graph Directed Markov Systems and Conformal Iterated Function 
Systems: Short Preliminaries}\label{CGDMS_Preliminaries}

This subsection has a preparatory character. It is needed for us in order to be able to formulate and to prove the main result, Theorem~\ref{t:thin}, of the next subsection. It establishes the Full Thin Annuli Property for essentially all conformal countable alphabet Iterated Function Systems. These iterated function systems will show up in later (applications, examples) sections too. There we will need them as our tool to prove the Full Thin Annuli Property for many other conformal dynamical systems such as conformal expanding repellers, rational functions of the Riemann sphere, and large subclasses of transcendental meromorphic functions from the complex plane to the Riemann sphere. An intermediate convenient tool, also interesting on its own, is that of (conformal) graph directed Markov systems introduced and systematically studied in \cite{GDMS}. These are considerable but quite natural generalizations of (conformal) countable alphabet iterated function systems.

Passing to strictly mathematical terms, let us define a graph directed Markov system (abbr. GDMS) relative to a directed multigraph $(V,E,i,t)$ and an incidence matrix $A$. As was indicated above, such systems were introduced and studied at length in \cite{GDMS}. A \emph{directed multigraph} consists of 

\begin{itemize} 

\item A finite set $V$ of vertices, 

\item A countable (either finite or infinite) set $E$ of directed edges, 

\item A map $A:E\times E\to \{0,1\}$ called an \emph{incidence matrix} on $(V,E)$,\

\item Two functions $i,t:E\to V$, such that $A_{ab} = 1$ implies $t(b) = i(a)$. 

\end{itemize}

\sp\noindent Now suppose that in addition, we have a collection of nonempty compact metric spaces $\{X_v\}_{v\in V}$ and a number $\lambda\in (0,1)$, and that for every $e\in E$, we have a one-to-one contraction $\phi_e:X_{t(e)}\to X_{i(e)}$ with Lipschitz constant $\le \lambda<1$. Then the collection
\[
\cS = \{\phi_e:X_{t(e)}\to X_{i(e)}\}_{e\in E}
\]
is called a \emph{graph directed Markov system} (or \emph{GDMS}). We now describe the limit set of the system $\cS$. For every $n \in \N$ let
\[
E_A^n:=\{\omega\in E^n:\forall (1\le j\le n-1) \; A_{\omega_j\omega_{j+1}}=1 \},
\]
and let $E_A^0$ be the set consisting of the empty word. Then let
\begin{equation}\label{120181110}
E_A^*:=\bigcup_{n=0}^\infty E_A^n
\end{equation}
and
\[
E_A^\infty:=\{\omega\in E^\N: \text{ every finite subword of } \omega   \text{ is in }   E_A^*\}.
\]
We use the commonly accepted symbol $\sg$ for the shift map
on $E^\N$. More precisely, $\sg:E^\N\to E^\N$ is defined by the formula:
$$
\sg\((\om_n)_{n=1}^\infty)\)=(\om_{n+1})_{n=1}^\infty.
$$
The space $E_A^\infty$ is forward invariant with respect to the shift map $\sg:E^\N\to E^\N$, i.e. $\sg\(E_A^\infty\)\subset E_A^\infty$ and we frequently consider also the dynamical system
$$
\sg:E_A^\infty\longrightarrow E_A^\infty.
$$
The union defining $E_A^*$ in the formula \eqref{120181110} is disjoint and for every $\omega\in E_A^*$ we denote by $|\omega|$ the unique integer $n$ such that $\omega\in E_A^n$; we call $|\omega|$ the \emph{length} of $\omega$. For each $\omega\in E_A^\infty$ and $n \in \N$, we write
\[
\omega|_n:=\omega_1\omega_2\ldots\omega_n\in E_A^n.
\]
For each $n \geq 1$ and $\omega \in E_A^n$, we let $i(\omega) = i(\omega_1)$ and $t(\omega) = t(\omega_n)$, and we let
\[
\phi_\omega:=\phi_{\omega_1}\circ\cdots\circ\phi_{\omega_n}:X_{t(\omega)}\to X_{i(\omega)}.
\]
For $\omega \in E^\infty_A$, the sets $\{\phi_{\omega|_n}\left(X_{t(\omega_n)}\right)\}_{n \geq 1}$ form a descending sequence of nonempty compact sets and therefore $\bigcap_{n \geq 1}\phi_{\omega|_n}\left(X_{t(\omega_n)}\right)\ne\emptyset$. Since for every $n \geq 1$, 
$$
\diam\left(\phi_{\omega|_n}\left(X_{t(\omega_n)}\right)\right)\le \lambda^n\diam\left(X_{t(\omega_n)}\right)\le \lambda^n\max\{\diam(X_v):v\in V\},
$$
we conclude that the intersection 
\[
\bigcap_{n \in \N}\phi_{\omega|_n}\left(X_{t(\omega_n)}\right)
\]
is a singleton and we denote its only element by $\pi(\omega)$. In this way we have defined a map
\[
\pi:E^\infty_A\longrightarrow \coprod_{v\in V}X_v,
\]
where $\coprod_{v\in V} X_v$ is the disjoint union of the compact sets $X_v \; (v\in V)$. The map $\pi$ is called the \emph{coding map}, and the set
\[
J = J_\cS = \pi(E^\infty_A)
\]
is called the \emph{limit set} of the GDMS $\cS$. The sets
\[
J_v = \pi(\{\omega \in E_A^\infty : i(\omega_1) = v\}) \;\; (v\in V)
\]
are called the \emph{local limit sets} of $\cS$.

\

We call the GDMS $\cS$ \emph{finite} if the alphabet $E$ is finite. Furthermore, we call $\cS$ \emph{maximal} if for all $a,b\in E$, we have $A_{ab}=1$ if and only if $t(b)=i(a)$.  In \cite{GDMS}, a maximal GDMS was called a \emph{graph directed system} (abbr. GDS).
Finally, we call a maximal GDMS $\cS$ an \emph{iterated function system} (or \emph{IFS}) if $V$, the set of vertices of $\cS$, is a singleton. Equivalently, a GDMS is an IFS if and only if the set of vertices of $\cS$ is a singleton and all entries of the incidence matrix $A$ are equal to $1$.

\begin{dfn}\label{definitionsymbolirred}
We call the GDMS $\cS$ and its incidence matrix $A$ \emph{finitely (symbolically) irreducible} if there exists a finite set $\Lambda\subset E_A^*$ such that for all $a, b\in E$ there exists a word $\omega\in\Lambda$ such that the concatenation $a\omega b$ is in $E_A^*$. $\cS$ and $A$ are called \emph{finitely primitive} if the set $\Lambda$ may be chosen to consist of words all having the same length. Note that all IFSs are finitely primitive.
\end{dfn}

Intending to pass to geometry and following \cite{GDMS}, we call a GDMS \emph{conformal} if for some $d\in\N$, the following conditions are satisfied:

\begin{itemize}
\item[(a)] For every vertex $v\in V$, $X_v$ is a compact connected subset of $\R^d$, and $X_v=\overline{\Int(X_v)}$.
\item[(b)] There exists a family of open connected sets $W_v \subset X_v \;\; (v\in V)$ such that for every $e\in E$, the map $\phi_e$ extends to a $C^1$ conformal diffeomorphism from $W_{t(e)}$ into $W_{i(e)}$ with Lipschitz constant bounded above by some number $\lambda<1$.

\item[(c)] There are two constants $L\ge 1$ and $\alpha>0$ such that for every $e\in E$ and every pair of points $x,y\in X_{t(e)}$,
\[
\left|\frac{|\phi_e'(y)|}{|\phi_e'(x)|}-1 \right| 
\le L\|y-x\|^\alpha,
\]
where $|\phi_\omega'(x)|$ denotes the scaling of the derivative, which is a linear similarity map.
\end{itemize}

\begin{rem}\label{p1.033101}
If $d\ge 2$ and a family $\cS = \{\phi_e\}_{e\in E}$ satisfies the conditions (a) and (b), then it also satisfies condition (c) with $\alpha=1$. When $d = 2$ this is due to the well--known Koebe's distortion theorem.
 When $d \geq 3$ it is due to \cite{GDMS} depending heavily on Liouville's representation theorem for conformal mappings, see \cite{IwaniecMartin} for a detailed development leading up to the strongest current version. If $d=1$, condition (c) is only automatic if the alphabet $E$ is finite and all contractions $\phi_e$ are of $C^{1+\epsilon}$ class. 
\end{rem}
\begin{rem}
We do emphasize that, unlike to \cite{GDMS}, in the above definition, and in all the results of this section, we do not need and we do not require any separation condition whatsoever. In particular even its weakest form
\[
\phi_a\left(\Int(X)\right)\cap \phi_b\left(\Int(X)\right)=\emptyset.
\]
for all $a,b\in E$ such that $a\ne b$, known as the Open Set Condition, is not assumed to hold. We also do emphasize that we do not impose any form of boundary regularity, in particular no Cone Condition of \cite{GDMS}.
\end{rem}

\subsection{Full Thin Annuli Property holds for (essentially all) Conformal IFSs}\label{sec:IFS}

In this subsection we establish the Full Thin Annuli Property for a large class of conformal countable alphabet IFSs.
\[
\cS = \{\phi_e:X\to X\}_{e\in E},
\]
where $X\subset\R^d$, $d\ge 1$. 

\begin{dfn}\label{as:d}
We say that the system $\cS$ is \emph{geometrically irreducible} if the limit set $J_\cS$ is not contained in any proper, i.e. of dimension $\le d-1$, real analytic sub-manifold; precisely: is not contained in a conformal image of any affine hyperspace or geometric round sphere of dimension $\le d-1$.
\end{dfn}
Throughout this whole Subsection~\ref{sec:IFS} we assume that the system $\cS$ is geometrically irreducible. For the sake of brevity we denote 
$$
D(\omega):=\diam(\phi_\omega(X))
$$
for all $\om\in E^*$. The Bounded Distortion Property tells us that
\begin{equation}\label{1_2016-01_28}
Q^{-1}D(\om)D(\tau)\le D(\om\tau)\le QD(\om)D(\tau) 
\end{equation}
for all $\om, \tau\in E^*$ and some constant $Q\ge 1$. In this section we consider a (really large) class, called $\mathcal M_E$, of Borel probability measures $\mu$ on the symbol space $E^\N$, determined by the following two requirements:
\begin{enumerate}[(A)]
\item Weak Independence: \label{as:weakind}
$$
P^{-1}\mu([\omega]) \mu([\tau]) 
\leq \mu([\omega\tau]) 
\leq P\mu([\omega])\mu([\tau])
$$
for some constant $P\ge 1$ and all $\om, \tau\in E^*$.
\item There exists $\beta >0$ such that $\displaystyle \sum_{e\in E} \frac{\mu([e])}{\diam^\beta(\phi_e(X))}<+\infty$. \label{as:c}
\end{enumerate}
\begin{rem}
All Gibbs measures, on the symbol space $E^\N$, introduced and considered in \cite{MU_Israel} are weakly independent, i.e. enjoy the property (A). It is easy to have abundance of measures satisfying the property (B); among them are the Gibbs states of all (geometrically most significant) potentials $E^\N\ni\om\longmapsto t\log\left|\phi_{\om_1}'(\pi(\sg(\om)))\right|\in\R$, where $t\ge 0$ is sufficiently large. 
\end{rem}
The main result of this subsection follows.
\begin{thm}\label{t:thin}
If $\cS =\{\phi_e\colon X\to X\}_{e\in E}$ is a geometrically irreducible conformal IFS, then for every $\mu\in\mathcal M_E$ the measure $\mu\circ\pi^{-1}$ on $J_\cS$ has the Thin Annuli Property with $\kappa=3$ (in fact this is true for any $\kappa>1$).
In other words:
\begin{align*}
\lim_{r\to 0} \frac{\mu\circ\pi^{-1}\left(B(x,r+r^3)\setminus B(x,r)
\right)}{\mu\circ\pi^{-1}\left(B(x,r)\right)}=0\mbox{\hskip1cm $\mu\circ\pi^{-1}$--a.e.}	
\end{align*}
\end{thm}
In order to ease notation, let us denote by $R(x,r,r^3)$ the annulus centred at $x$ with inner radius $r>0$ and outer radius $r+r^3$, i.e. 
\[
R(x,r,r^3):=B(x,r+r^3)\setminus B(x,r). 
\]
The proof of Theorem~\ref{t:thin} consists of several steps listed below, and it has been strongly    influenced by the techniques of \cite{DFSU}. 
For the sake of brevity we denote $\hat\mu:=\mu\circ\pi^{-1}$.
\begin{lem}\label{l:cyl}
There exist constants $\rho>0$, $H<\infty$ and a finite set $F\subset E^*$ such that for any $x\in J_\cS$, any radius $0<r<\rho$, and any finite word $\omega\in E^*$, with diameter $D(\omega)\geq Hr^3$, there exists a word $\tau\in F$ such that $\pi([\omega\tau])$ does not intersect 
$R(x,r,r^3)$. In symbols:
\[
\pi([\omega\tau])\cap R(x,r,r^3)=\emptyset.\]
\end{lem}
\begin{lem} \label{l:est}
There exist $\alpha>0$, $C<\infty$, $\rho>0$ such that for all $0<r<\rho$, $x\in J_\cS$ and any finite word $\omega\in E^*$, with diameter $D(\omega)\geq r^2$, we have
\begin{equation}
\mu\left([\omega]\cap \pi^{-1}(R(x,r,r^3))\right)\leq C r^\alpha\mu([\omega]).
\end{equation}
\end{lem}
\begin{lem}\label{l:T}
For any numbers $0<A<B$ define the set
\begin{equation}
T_{A}^B:=\left\{\omega\in E^\N\colon \forall_{k\in \N} \, D(\omega|_k)\notin (A,B)\right\}.
\end{equation}
Then there exists $C<\infty$ for which $\mu(T_{A}^B) \leq C \left(\frac{A}{B}\right)^\beta \ln\left(\frac{\diam X}{A}\right)$, where $\beta$ is the constant from condition \eqref{as:c} from the definition of the space $\mathcal{M}_E$.
\end{lem}
\begin{lem}\label{l:S}
Let $\nu$ be an arbitrary Borel probability measure defined on some bounded Borel set $X\subset \R^d$. Let $F$ be a measurable subset of $X$. Define
\begin{equation}
S(F,c,\rho):=\left\{x\in X\colon \nu(B(x,\rho)\cap F) > c \nu(B(x,\rho))\nu(F)\right\}.
\end{equation}
Then for any numbers $c,\rho>0$ we have $\nu(S(F,c,\rho))\leq M/c$, where $M$ is some constant depending only on the space $X$. 
\end{lem}
\begin{proof}[Proof of Lemma~\ref{l:cyl}] 
Assume without loss of generality that $E=\N$. Seeking a contradiction suppose that there exist a sequence $(r_n)_{n=1}^\infty\searrow 0$, a sequence $x_n\in J_\cS$, $n\in\N$, and a sequence of finite words $\omega^{(n)}\in E^*$ with diameters satisfying
\begin{equation}\label{1_2015_12_04}
D(\omega^{(n)})\geq nr_n^3
\end{equation}
such that for every $\tau\in \{1,2,\ldots,n\}^n$ the \emph{cylinder} $\pi([\omega^{(n)}\tau])$ intersects $R(x_n,r_n,r_n^3)$. Let us denote 
$$
R_n:=R(x_n,r_n,r_n^3) \quad \mbox{and} \quad  S_n:=\bd B(x_n,r_n)=\{x\in\R^d:||x-x_n||=r_n\}.
$$
Take then any sequence of similarities $T_n$, $n\ge 1$, for which $0\in T_n(\pi([\omega^{(n)}]))$ and $|T_n'|=(D(\omega^{(n)}))^{-1}$ for all $n\ge 1$. Note that $(T_n\circ \phi_{\omega^{(n)}})_{n=1}^\infty$ is a bounded equicontinuous sequence of conformal maps with derivatives uniformly bounded from above and uniformly separated from zero. Actually, in dimension $d=1$ the conformality is not needed, making the proof easier. 

Applying Ascoli-Arzela Theorem and passing to an appropriate subsequence we will have that the sequence $(T_n\circ \phi_{\omega^{(n)}})_{n=1}^\infty$ converges uniformly on $X$ to a conformal map $U:X\to\R^d$. Now, working with the one-point (Alexandrov) compactification $\hat\R^d$ of $\R^d$, with $\infty$ as the compactifying point, endowing $\hat\R^d$ with spherical metric, and then the collection $\cK_d$ of non-empty compact subsets of $\hat\R^d$ with the corresponding Hausdorff metric $d_H$, we see that the collection $\Ga$ of all geometric spheres of $\hat\R^d$, including the spheres containing infinity (hyperplanes) and singletons, forms a compact subset of $\cK_d$. Since $T_n(S_n)\in\Ga$, passing to a subsequence, we can therefore assume without loss of generality that $T_n(S_n)$ converges in the Hausdorff metric $d_H$ to some element $Q\in\Ga$. 

Depending on actual sizes of $D(\omega^{n)})$, the limit object $Q$ may be either a sphere -- the case if $D(\omega^{(n)})
\asymp r_n$,  a point in $\mathbb R^d$ -- the case if $D(\omega^{(n)})/r_n\nolinebreak\to\nolinebreak \infty$, or a line in $\mathbb R^d$ -- which is so if $D(\omega^{(n)})/r_n\to 0$.  In all three cases the ratio of the outer and inner radii of the annulus  $T_n(R_n)$ converges to one, as $\frac{r+r^3}{r}\to 1$ when $r\to 0$. 

In the first two cases, this immediately proves that also $\lim_{n\to\infty}T_n(R_n)=Q$. In the third case, we  need to use additionally (\ref{1_2015_12_04}), to conclude that both spheres bounding the annulus $R(x_n, r_n, r_n+r_n^3)$, after rescaling   by $(D(\omega^{n)}))^{-1}\le\frac{1}{nr_n^3}$   tend to the same line in $\mathbb R^d$.

So, finally, in all three cases we may conclude that 
\begin{equation}\label{1_2015_12_08}
\lim_{n\to\infty}T_n(R_n)=Q.
\end{equation}
Observe also that by Definition \ref{as:d} for every $M\in \Ga$ there exists a point $w_M\in J_\cS$ such that $\dist(w_M,M)>0$. Writing the $w_{U^{-1}(Q)}=\pi(\xi)\in J_\cS$, $\xi\in E^{\N}$, we have that
$$
\dist(\pi(\xi),U^{-1}(Q))>0.
$$
We therefore conclude that there exists $k\ge 1$ such that 
\begin{equation}\label{2_2015_12_08}
\dist\big(\pi([\xi|_k]),U^{-1}(Q)\big)>0.
\end{equation}
Consider now only integers $n\ge k$ so large that all letters forming $\xi|_{k}$ belong to $\{1,2,\ldots,n\}$. By our contrary hypothesis
$$
\phi_{\omega^{(n)}}\big(\phi_{\xi|_k}(J_\cS)\big)\cap R_n=\phi_{\omega^{(n)}\xi|_k}(J_\cS)\cap R_n\ne\emptyset.
$$
Fix an arbitrary $z_n\in J_\cS$ such that $\phi_{\omega^{(n)}\xi|_k}(z_n)\in R_n$. Passing to a subsequence we may assume without loss of generality that $\lim_{n\to\infty}z_n=z\in X$ for some point $z\in \overline J_\cS$. Then, invoking also \eqref{1_2015_12_08}, we get that
$$
\lim_{n\to\infty}T_n\circ \phi_{\omega^{(n)}}(\phi_{\xi|_k}(z))
=U(\phi_{\xi|_k}(z))\in Q.
$$
Hence $\phi_{\xi|_k}(z)\in U^{-1}(Q)$, and as $\phi_{\xi|_k}(z)\in\overline{\pi([\xi|_k])}$, this contradicts \eqref{2_2015_12_08} and finishes the proof of our lemma.
\end{proof}
%
\begin{proof}[Proof of Lemma \ref{l:est}]
Take $\rho$, $H$ and the set $F$ given by Lemma \ref{l:cyl}. First of all, observe that, because the lengths of all words in $F$ are uniformly bounded above, by taking an iterate of the system $\cS$ we may assume that $F\subset E$ (instead of $E^*$). Fix $x\in J_\cS$, $0<r<\rho$, and denote $R:=R(x,r,r^3)$.

We will in fact prove a stronger fact; namely that with no restrictions on $D(\omega)$
\begin{equation}
\mu\left([\omega]\cap \pi^{-1}(R)\right) \leq\left(\frac{Hr^3}{D(\omega)}\right)^\alpha\mu([\omega]),
\label{eq:ic}
\end{equation}
for all $\omega\in E^*$. This will trivially prove the lemma as its hypotheses require that $D(\omega)\geq r^2$. So, we now focus on the proof of (\ref{eq:ic}). First note that if $D(\omega)\leq Hr^3$, then inequality (\ref{eq:ic}) is trivial. Also for all $n\ge 1$ big enough and all $\omega\in E^n$ we have $D(\omega)\leq Hr^3$.
Now let us work \emph{from the bottom upwards}. 
Take a cylinder $[\omega]$ such that (\ref{eq:ic}) is already proven for all subcylinders $[\omega e]$, $e\in E$. We have
\[
\mu([\omega]\cap \pi^{-1}(R))
	= \sum_{e\in E} \mu([\omega e]\cap \pi^{-1}(R)).
\]
Applying Lemma~\ref{l:cyl}, we may drop at least one element of this sum, say $b\in F$, to get
\begin{align*}
\mu([\omega]\cap \pi^{-1}(R))
&= \sum_{E\ni a\neq b}\mu([\omega a]\cap \pi^{-1}(R))  \\
&  \leq \sum_{E\ni a\neq b} \left(\frac{Hr^3}{D(\omega a)}\right)^\alpha\mu([\omega a])  \\
&=\left(Hr^3\right)^\alpha \sum_{E\ni a\neq b} \frac{\mu([\omega a])}{(D(\omega a))^\alpha},
\end{align*}
where we used the estimate (\ref{eq:ic}) for every cylinder $[\omega a]$. In order to prove the required  inequality we need to have
\[
\sum_{E\ni a\neq b} \frac{\mu([\omega a])}{(D(\omega a))^\alpha} 
\leq \frac{\mu([\omega])}{(D(\omega))^\alpha} 
= \sum_{a\in E} \frac{\mu([\omega a])}{(D(\omega))^\alpha},
\]
where the equality sign trivially holds. Simplifying this gives 
\[
\sum_{E\ni a\neq b} \left(\left(\frac{D(\omega)}{D(\omega a)}\right)^\alpha-1\right) \mu([\omega a])  
\leq \mu([\omega b]).
\]
Applying Bounded Distortion Property \eqref{1_2016-01_28} and Weak Independence of $\mu$, i.e. condition (A), we see that it is thus enough to prove that
\[
\sum_{E\ni a\neq b} \left(\left(\frac{QD(\omega)}{D(\omega)D(a)}\right)^\alpha-1\right) P\mu([\omega])\mu([a])
\leq P^{-1} \mu([\omega])\mu([b]).
\]
Recall that $b$ was chosen from a finite set so $P^{-2}\mu([b])$ is bounded away from zero, say $P^{-2}\mu([b])>\delta$ for some fixed $\delta>0$. Simplifying again, we see that it is enough to prove
\[
\sum_{E\ni a\neq b} \left(\left(\frac{Q}{D(a)}\right)^\alpha-1\right)\mu([a])\leq \delta. \]
Therefore, it is enough to have
$$
\sum_{E\ni a\neq b}\frac{\mu([a])}{D(a)^\alpha}
\leq Q^{-\alpha}\delta+Q^{-\alpha}\sum_{E\ni a\neq b}\mu([a]).
$$
But since, by Assumption \eqref{as:c}, the series on the left-hand side of this formula converges for all $\alpha>0$ small enough. Its sum tends to $\sum_{E\ni a\neq b}\mu([a])$ as $\alpha\to 0$ and using a dominated convergence theorem we get that this formula will hold for all $\alpha>0$ small enough. Thus the proof is complete.
\end{proof}

\begin{proof}[Proof of Lemma \ref{l:T}] 
First, divide $T_A^B$ into disjoint subsets (for $k=0,1\ldots$)  
\[
T_A^B(k)=\left\{\omega\in E^\N\colon D(\omega|_{k+1})\leq A<B\leq D(\omega|_k)\right\}.
\]
For any cylinder $D(\omega|_k)\leq \lambda^k\diam(X)$, so for any $n\geq N:=\log_\lambda\left(\frac{A}{\diam X}\right)$ we have $D(\omega|_n)\leq A$ and $T_A^B(n)=\emptyset$. This allows us to write
\begin{equation}
\mu(T_A^B)\leq \sum_{n=0}^N \mu\left(T_A^B(n)\right).
\label{eq:l41}
\end{equation}
Now, fix $0\le k\le N$ and $\omega\in E^\N$. If $D(\omega|_k)<B$, then $\mu\left(T_A^B(k)\cap[\omega|_k]\right)=0$. If $D(\omega|_k)\geq B$, then
\[
\mu\left(T_A^B(k)\cap[\omega|_k]\right)=\sum_e\mu([\omega|_ke]),
\]
where the sum is taken over those $e\in E$ for which $D(\omega|_ke)\leq A$. Applying the Weak Independence of $\mu$, i.e. condition (A) and Bounded Distortion \eqref{1_2016-01_28}, we further get
\[
\mu\left(T_A^B(k)\cap[\omega|_k]\right)
\leq \!\sum_{D(\omega|_ka)\leq A}\! P\mu([\omega|_k])\mu([a])
\leq \!\!\!\sum_{Q^{-1}D(\omega|_k)D(a)\leq A}\!\!\!P\mu([\omega|_k])\mu([a]),
\]
and using the fact that $D(\omega|_k)\geq B$, this gives
\begin{equation}\label{eq:l42}
\mu(T_A^B(k)\cap[\omega|_k])
\leq P\mu([\omega|_k])\sum_{D(a)\leq QA/B}\mu([a]).
\end{equation}
By Assumption~\eqref{as:c} we may write:
\begin{align*}
+\infty>Z
&:=\sum_{a\in E} \frac{\mu([a])}{D(a)^\beta}
 \geq \sum_{D(a)\leq QA/B} \frac{\mu([a])}{D(a)^\beta}
 \geq \sum_{D(a)\leq QA/B} \frac{\mu([a])}{(QA/B)^\beta} \\
&= \sum_{D(a)\leq QA/B} \mu([a])\left(\frac{B}{QA}\right)^\beta.
\end{align*}
Combining this estimate with (\ref{eq:l42}) gives
\[
\mu\left(T_A^B(k)\cap[\omega|_k]\right)
\leq P\mu([\omega|_k]) \cdot Z\left(\frac{QA}{B}\right)^\beta,
\]
and summing over all cylinders $[\omega|_k]$, this gives $\mu(T_A^B(k))\leq C(A/B)^\beta$ with some constant $C$. Finally applying (\ref{eq:l41}), we get 
\[\mu(T_A^B)\leq \log_\lambda\left(\frac{A}{\diam X}\right) \cdot C\left(\frac{A}{B}\right)^\beta.\qedhere\]
\end{proof}
\begin{proof}[Proof of Lemma \ref{l:S}]
Set $S:=S(F,c,\rho)$. By Besicovitch's Covering Theorem there exists a covering of $S$ with balls $B(x_i,\rho)$, $i\in I$, all centred at $S$, with finite multiplicity $M_d$ depending only on the dimension $d$. The following estimate uses first, the definition of $S$ and then the multiplicity, bounded by $M_d$, of the covering. 
\[
\nu(S)\leq \sum_{i\in I} \nu(B(x_i,\rho)) \leq \sum_{i\in I} \frac{\nu(B(x_i,\rho)\cap F)}{c\nu(F)}\leq \frac{M_d\nu(F)}{c\nu(F)} = \frac{M_d}{c}.\qedhere
\]
\end{proof}
In the final proof of this section we also use Proposition \ref{prop2:BS} (proved in \cite{BS}). Recall that it is an immediate consequence of, much stronger, Theorem~\ref{thm:BSs}.
\begin{proof}[Proof of Theorem \ref{t:thin}] 
We will show that for $\hat\mu$ almost every $x\in J_\cS$ and all sufficiently small radii $r>0$ we have that, for some $\gamma>0$,
$$
\hat\mu(R(x,r,r^3))\leq C\hat\mu(B(x,r))r^\gamma.
$$ 
First, using notation from Lemmas \ref{l:T} and \ref{l:S} define 
$T_n:=T_{4^{-n}}^{2^{-n}}$, $n\ge 1$ and denote $\hat{T_n}=\pi(T_n)$. Put
\[
S_n:=S(\hat{T}_n,n^2, 4\cdot 2^{-n}).
\] 
Lemma~\ref{l:S} gives that $\hat\mu(S_n)\leq M/n^2$ and so $\sum_n \hat\mu(S_n)<\infty$. Thus the Borel--Cantelli Lemma applies to tell us that for  $\hat\mu$ almost every $x\in J_\cS$ there exists an integer $K(x)\ge 1$ such that $x \notin S_k$ for all $k\geq K(x)$.
Fix $x\in J_\cS$ with such property, i.e. an arbitrary $x$ produced by the Borel--Cantelli Lemma. For any $n\ge 1$ define the set 
\begin{equation}
C_n=\left\{[\omega]\in E^*\colon D(\omega)\leq 2^{-n}<D(\omega|_{|\omega|-1})\right\}.
\end{equation}  
Now, take any $0<r\leq 2^{-(K(x)+1)}$. Define $n\ge 1$ so as to satisfy the inequalities $2^{-n-1} < r \leq 2^{-n}$. Then
\begin{equation}\label{1_2015_12_10}
n\ge K(x).
\end{equation}
Denote the annulus $R(x,r,r^3)$ by $R$ and cover $\pi^{-1}(R)$  by cylinders from $C_n$. We estimate the measure 
$$
\begin{aligned}
\hat\mu(R)
&=\mu\circ\pi^{-1}(R)
 \leq \sum_{[\omega]}\!{}^* \mu([\omega]\cap \pi^{-1}(R)) \\
&\leq \underbrace{\sum_{[\omega]\subset T_n}\!\!\!\!{}^* \  \mu([\omega]\cap \pi^{-1}(R))}_{I}+\!\!\!\underbrace{\sum_{[\omega]\cap T_n=\emptyset}
 \!\!\!\!\!\!{}^* \ \mu([\omega]\cap \pi^{-1}(R))}_{II},
\end{aligned}
$$
where the ${}^*$ indicates that the corresponding sum above is taken over all cylinders $[\omega]\in C_n$ intersecting $\pi^{-1}(B(x,r+r^3))$.
Recall that for such cylinders $D(\omega)<2r$, and as $r+r^3\leq 2r$, the cylinder $[\omega]$ is contained in the set $\pi^{-1}(B(x,4r))$. So 
$$
I\leq \sum_{[\omega]\subset T_n}\!\!\!\!{}^* \   \mu([\omega]) 
\leq \mu\left(T_n \cap \pi^{-1}(B(x,4r))\right) 
\leq \mu\left(T_n \cap \pi^{-1}(B(x,4\cdot 2^{-n}))\right).
$$
Now, first straightforward from the definition of $S_n$, and from the fact that, because of \eqref{1_2015_12_10}, $x\notin S_n$, then by applying Lemma~\ref{l:T}, we get that
$$
\begin{aligned}
I
&\leq n^2\mu\left(\pi^{-1}(B(x,4\cdot 2^{-n}))\right)\mu(T_n) \\
&\leq n^2\hat\mu\left(B(x,4\cdot 2^{-n})\right) C \left(\frac{4^{-n}}{2^{-n}}\right)^\beta \ln\Big(\frac{\diam X}{4^{-n}}\Big)\\
&\leq n^2\hat\mu\left(B(x,4\cdot 2^{-n})\right)\cdot \widehat{C} n 2^{-n\beta}\\
&\leq \widetilde{C}\hat\mu\left(B(x,8r)\right)r^{\beta/2}
\end{aligned}
$$
with appropriate constants $\widehat{C}$ and $\widetilde{C}$. 

Finally, we apply the estimate of Proposition~\ref{prop2:BS} with $\varepsilon=\beta/4$ to get
$$
I\leq \widetilde{C}\hat\mu\left(B(x,r)\right)r^{-\varepsilon}r^{\beta/2} \leq \widetilde{C}\hat\mu\left(B(x,r)\right)r^{\beta/4} 
$$
which completes the estimate of the first sum, i.e. the one labelled by $I$.

Now, observe that if $[\omega]\cap T_n=\emptyset$, then $D(\omega)\geq 4^{-n}\geq r^2$ and so we may first apply Lemma~\ref{l:est}, and then Proposition~\ref{prop2:BS} with $\varepsilon=\alpha/2$ to estimate as follows:
$$
\begin{aligned}
II
&\leq \sum_{[\omega]\cap T_n=\emptyset}\!\!\!\!\!\!{}^* \  Cr^\alpha\mu([\omega])
\leq Cr^\alpha \hat\mu(B(x,4r)) \\
&\leq Cr^\alpha \hat\mu(B(x,r))r^{-\varepsilon} \leq  Cr^{\alpha/2} \hat\mu(B(x,r)).
\end{aligned}
$$
This completes the upper estimate of $II$ and finishes the entire proof. 
\end{proof}
\section{Applications and Examples: Exponential One Laws}\label{sec:ex}
\subsection{Expanding Repellers}\label{DEMER}
In this subsection we deal with the class of, not assumed to be conformal, expanding repellers. The main result of this short subsection is Theorem~\ref{thm:ER}. It was proved  in  \cite{saussolphys} (see also  \cite{R} for its random counterpart) with the extra hypothesis that a sufficiently strong version of the thin annuli property holds. The main point in our approach is that we do not assume any form of thin annuli property. We proved it: this is Theorem~\ref{Main_Thm_3}. Because of the aforementioned papers, having Theorem~\ref{Main_Thm_3}, we could have actually skipped the actual proof of Theorem~\ref{thm:ER}, merely referring to them. We however provide it as a prelude to more technically involved further sections, for the sake of completeness, convenience of the reader, and because this proof is quite short. 
 
In what follows 
we will need the classical concepts of topological pressure, variational principle, and equilibrium states. We bring them up now. Let $X$ be a compact metrizable space, $T\colon X\to X$ be a continuous map, and $\varphi \colon X\to\mathbb{R}$ be a 
continuous function. We denote by $\P(\varphi)$ its topological pressure with respect to the dynamical system given by $T$, see for example \cite{PUbook} for the definition and properties.
One of the most important of these properties is the following formula, commonly referred to as the Variational Principle.
\begin{equation}\label{VP}
\P(\varphi)=\sup\Big\{\h_\mu(T)+\int_X\varphi \,d\mu \Big\},
\end{equation}
where the supremum is taken over all Borel probability $T$--invariant measures on $X$. Any measure for which the  supremum is attained is called an equilibrium state of $\varphi$.

We now provide the definition of expanding repellers.

\begin{dfn}\label{exprep2}
Let $U$ be an open subset of $\R^d$, $d\ge 1$. Let $J$ be a compact subset of $U$. Let $T:U\to\R^d$ be a $C^{1+\epsilon}$--differentiable map. The map $T$ is called an expanding repeller, if the following conditions are satisfied:
\begin{enumerate}
\item{}$T(J)=J$,
\item
for every $z\in J$ the derivative $T'(z):\R^d\to\R^d$ is invertible and the norm of its inverse is smaller than $1$.
\item there exists an open set  $V$ such that $\overline{V}\subset U$ and
$$
J=\bigcap_{k=0}^\infty T^{-n}(V).
$$
\item the map $T|_{J}:J\to J$ is topologically transitive.
\end{enumerate}
Note that $T$ is not required to be one-to-one; in fact usually it is not. Abusing  notation slightly we refer to the set $J$ alone as an expanding repeller. In order to use the uniform terminology we also call $J$ the limit set of $T$. 
\end{dfn}

One of the basic concepts associated with expanding repellers is this. 
\begin{dfn}\label{3.5.1}
A finite cover $\cR=\{R_1,\ldots,R_q\}$ of $X$ is
said  to be a Markov partition of the space $X$ for the mapping $T$ if
the following conditions hold. 
\begin{itemize}
\item[(a)] \ $R_i=\overline{\Int R_i}$ \ for all $i=1,2,\ldots,q$.
\item[(b)] \ $\Int R_i\cap \Int R_j = \es$ \ for all $i\ne j$.
\item[(c)] \ $\Int R_j\cap T(\Int R_i) \ne \es \  \Longrightarrow \  R_j\subset T(R_i)$ \ for all $i,j=1,2,\ldots,q$.
\end{itemize}
\end{dfn}
\noindent The elements of a Markov partition will be called cells in the sequel. The existence of Markov partitions is guaranteed by the following theorem whose proof can be found for instance in \cite{PUbook}. 
\begin{thm}\label{t_Markov_Partitions}
Any expanding repeller $T\colon J\to J$ admits Markov partitions of arbitrarily small diameters.
\end{thm}

Another crucial theorem about expanding repellers is the following, see \cite{PUbook} for a proof. 

\begin{thm}\label{t120181109}
If $T\colon J\to J$ be an expanding repeller and $\psi\colon J\to\R$ is a H\"older continuous potential, then there exists $\mu_\psi$, a unique equilibrium state for $\psi$ with respect to $T$. 
\end{thm}

The equilibrium state $\mu_\psi$ is also a unique Gibbs state of $\psi$. A definition of them can be again found in \cite{PUbook}; we will not need it here. We are now ready to state and prove the main result of this subsection. 

\begin{thm}\label{thm:ER}
Let $T\colon J\to J$ be an expanding repeller, let $\psi\colon J\to\R$ be a H\"older continuous potential, and let $\mu_\psi$ be the corresponding equilibrium (Gibbs) state. Then the measure--preserving dynamical system $\(T,\mu_\psi\)$ is Weakly Markov. In particular, the exponential one laws of Theorem~\ref{wykl_LM2} hold. 
\end{thm}
\begin{proof}
We shall check that the system satisfies the requirements of  Definition~\ref{dfn:Weakly_Markov} defining Weakly Markov systems. Property (i) of this definition for the dynamical system $\(T,\mu_\psi\)$ has been proved in \cite{PUbook}. Property (ii) also has been proved therein. For property (iii) we also use \cite{PUbook} and Markov partitions discussed above. 
We aim to show that these partitions fulfil the requirements of Definition~\ref{def:weak_partition}, i.e. the Weak Partition Existence Condition.

Towards this end fix $\delta>0$ so small that for every $x\in X$ and every $n\ge 0$ there exists $T_x^{-n}\colon B(T^n(x),4\delta)\to\R^d$, a unique continuous branch of $T^{-n}$ sending $T^n(x)$ to $x$. Theorem~\ref{t_Markov_Partitions} guarantees us the existence of 
$\cR=\{R_1,\ldots,R_q\}$, a Markov partition of $T$ with all cells of diameter smaller than $\delta$. It is not hard to see (proof in \cite{PUbook}), that any two distinct elements of $\cR$ intersect along a set of $\mu_\psi$ measure zero. So, we can treat $\cR$ as an ordinary partition. Its entropy $\h_{\mu_\psi}(T,\cR)$ is finite since the partition $\cR$ is finite, and this entropy is positive for all $\delta>0$ small enough since $\h_{\mu_\psi}(T)>0$. 

We now shall check that formula \eqref{borel_cantelli} holds. Fix one element 
$\xi\in R_1$.
Now fix $R>0$ so small that 
\begin{equation}\label{2_2016_09_08}
B(\xi,2R)\sbt R_1.
\end{equation}
Since $\mu_\psi\(B(\xi,R)\)>0$ (as $\mu_\psi$ has full topological support in $J$), it follows from ergodicity of $\mu_\psi$ and Birkhoff's Ergodic Theorem that for $\mu_\psi$--a.e. $z\in J$ there exists an infinite increasing sequence $(n_j)_{j=1}^\infty$ of integers such that
\[
T^{n_j}(z)\in B(\xi,R)
\mbox{\qquad for all $j\ge 1$ and \qquad}
\lim_{j\to\infty}\frac{n_{j+1}}{n_j}=1.
\]
So, there exists a constant $A\ge 1$ such that $n_{j+1}\le An_j$ for all $j\ge 1$. One consequence of such choice of $z$ is that
$$
z\notin\bu_{n=0}^\infty T^{-n}\bigg(\bu_{i=1}^q\partial\cR_i\bigg).
$$
In particular all elements $\cR^n(z)$, $n\ge 0$, are well-defined and $\cR$ being a Markov partition yields
\begin{equation}\label{1_2016_09_08}
\cR^n(z)=T_z^{-n}(\cR(T^n(z))).
\end{equation}
Fix an integer $k>n_1$. There exists a unique integer $j\ge 2$ such that 
$
n_{j-1}<k\le n_j
$. 
We then have $k>A^{-1}n_j$, and with $L\ge 1$ being a Lipschitz constant of $T$, looking up at \eqref{1_2016_09_08} and \eqref{2_2016_09_08}, we get that
\begin{equation}\label{eq:exprepest}
\begin{aligned}
\cR^k(z)
&\supset \cR^{n_j}(z)
=T_z^{-n_j}(\cR(T^{n_j}(z)))
=T_z^{-n_j}(R_1)\\
&\supset T_z^{-n_j}\(B(T^{n_j}(z),R)\) 
\supset B(z,L^{-n_j}R)\\
&\supset B(z,L^{-Ak}R) =B\(z,R\exp\(-A\log Lk\)\).
\end{aligned}
\end{equation}
Also $B\(z,R\exp\((-A\log L)k\)\)\supset B\(z,\exp\(-2A\log Lk\)\)$ for all $k\ge 1$ large enough. So,
\begin{equation}\label{1rv1}
\cR^k(z)
\supset B\(z,\exp\(-2A\log Lk\)\)
\end{equation}
for all such $k$, say $k\ge k(z)\ge n_1$. On the other hand, obviously there exists some $\chi^*(z)>0$ so large that
$$
\cR^k(z)
\supset B\(z,\exp(-\chi^*(z)k)\)
$$
for all $k=0,1,\ldots, k(z)-1$. In conjunction with \eqref{1rv1} this gives that
$$
\cR^k(z)
\supset B\Big(z,\exp\(-\max\big\{2A\log L,\chi^*(z)\big\}k\)\Big)
$$
for all $k\ge 0$, and formula \eqref{borel_cantelli} is verified. Thus, Weakly Markov Property is established. 
\end{proof}


\subsection{Equilibrium Measures (States) for Holomorphic Endomorphisms of Complex Projective Spaces}\label{projective_spaces}

Let $f:\mathbb P^k\to\mathbb P^k$ be a holomorphic endomorphism of a complex
projective space $\mathbb P^k$, $k\ge 1$, and let $J(f)$ be the Julia set of
$f$, which is commonly defined to be the topological support of the (unique) Borel probability $f$--invariant measure of maximal entropy. Generally, this system is not conformal, although sometimes it is, for example if $k=1$, the case dealt with in Section~\ref{holder_equil}. Let $\varphi:J(f)\to\R$ be a H\"older continuous function. It was proved in \cite{UZ_projective} that if  
$$
\sup(\varphi)-\inf(\varphi)<\kappa_f,
$$ 
where $0<\kappa_f\le \log d$ is some constant depending on the map $f$,  then $\varphi$ admits a
unique equilibrium state $\mu_\varphi$ on $J(f)$.  Further strong stochastic properties of the measure $\mu_\varphi$ were established in \cite{suz_projective}. A potential $\varphi$ satisfying the above condition is called admissible. 

Before proving Theorem~\ref{thm: projective_equil} below, the main result of this section, we formulate a technical, now rather standard result; see in particular \cite{suz} Proposition~10, for a similar statement.

\begin{prop}\label{prop:partition}
If $\mu$ is a finite Borel measure in $\mathbb P^k$, then for every $\delta>0$ there exists a finite partition $\alpha=\{U_i\}_{i\in I}$ of $\mathbb P^k$ with $\diam(U_i)<\delta$ (the diameter being calculated with respect to the Fubini--Study metric) for all $i\in I$, and such that 
\begin{equation}\label{boundary}
\mu\left(\bigcup_{j\in I} B(\partial U_j,r)\right)\le r^{1/2}
\end{equation}
for all sufficiently small $r>0$.
In fact, the  number $1/2$ in formula \eqref{boundary} can be replaced  by any positive number smaller than $1$.
\end{prop}

\fr For every $z\in\mathbb P^k$ denote by $\alpha(z)$ the only element of $\alpha$ containing $z$.
Now, we shall prove the following main result of this section.

\begin{thm}\label{thm: projective_equil}
Let $f:\mathbb P^k\to \mathbb P^k$, $k\ge 1$, be a holomorphic endomorphism of a complex projective space $\mathbb P^k$ of degree $d\geq 1$. Let $\varphi:J(f)\to\mathbb{R}$  be an admissible potential, and let $\mu_\varphi$ be its unique equilibrium state. Then  $(J(f), f, \mu_\varphi)$ forms a Weakly Markov system. Consequently, Theorem~\ref{wykl_LM2} holds for this system.
\end{thm}
\begin{proof} 
We shall check that the system satisfies the requirements of Definition~\ref{dfn:Weakly_Markov} defining Weakly Markov 
systems. Item (i) of this definition follows 
from Theorem~7.6 in \cite{suz_projective}. Item (ii), i.e. positive lower pointwise dimension, can be deduced from much more precise estimate of the lower pointwise dimension of some $f$--invariant Borel probability measures obtained in \cite{dupont}, Theorem A. Indeed, note that the hypothesis of this theorem, 
$$
\h_\mu(f)>(k-1)\log d,
$$
is fulfilled for our system since,
$$
\begin{aligned}
\h_{\mu_\varphi}(f)+\int_J\varphi\, d\mu_\varphi
&=\P(\varphi)\ge \h_m(f)+\int_J\varphi\, dm\ge k \log d+\inf(\varphi)\\
&>k \log d+\sup(\varphi)-\log d \\
&=\sup\varphi+(k-1)\log d,
\end{aligned}
$$ 
where we denoted by $m$ the measure of maximal entropy of $f$. Hence, 
$$
\h_{\mu_\varphi}(f)
>(k-1)\log d+\(\sup(\varphi)-\int\varphi \,d\mu_\varphi\)
\ge(k-1)\log d,
$$
as required.

By virtue of Remark~\ref{rem:part}, in order to prove property (iii) of Definition~\ref{dfn:Weakly_Markov}, it is enough to check that the Weak Partition Existence Condition holds.  We do it now.
Indeed, Proposition~\ref{prop:partition} provides  a partition $\alpha$ with elements of arbitrarily small diameter, satisfying the estimate ~\eqref{boundary}. If $\max\{\diam(U_i):i\in I\}$ is sufficiently small, then $h_{\mu_\varphi}(f,\alpha)>0$ as can be immediately seen by combining Shannon--Breiman--McMillan Theorem together with a local  entropy formula in  \cite{Katok_Brin}.

Now, we shall argue that condition \eqref{borel_cantelli} in the Definition ~\ref{def:weak_partition} (Weak Partition Existence Condition) is satisfied.  Fix $\beta>0$ arbitrary (later it will be needed to be sufficiently large) and for every integer $n\ge 1$  put 
$$
A_n:=f^{-n}\(\bigcup_{i\in I} B(\partial U_i,e^{-\beta n})\).
$$
Using the estimate \eqref{boundary} and $f$--invariance of measure $\mu_\varphi$, we see that 
$$
\mu_\varphi(A_n)\le e^{-(\beta/2) n}
$$
for every $n\ge 1$ provided that $\beta>0$ is sufficiently large.  Since the series $\sum_{n\ge 1}e^{-(\beta/2)n}$ converges, Borel--Cantelli Lemma thus applies and it tells us that 
for $\mu_\varphi$--a.e. $x\in J(f)$  there exists an integer $N=N(x)\ge 1$ such that for all integers $n\ge N$
\begin{equation}\label{1_2016_04_14}
B\(f^n(x),e^{-\beta n}\)\subset \alpha(f^n(x)).
\end{equation}
Keep such an $x$ and assume in addition that 
$$
x\notin \bigcup_{n=0}^\infty\bigcup_{i\in I}f^{-n}(\partial U_i).
$$
The set $A$ of such all points $x\in J(f)$ is of full measure, i.e. $\mu_\varphi(A)=1$. For all integers $n\ge N=N(x)$ denote by 
$\alpha^n_N(x)$ the only element of the partition $\bigvee _{k= N}^n f^{-k}(\alpha)$, containing the point $x$. Similarly, denote by  $\alpha^{N-1}_0(x)$ the only element of the partition $\bigvee _{k=0}^{N-1} f^{-k}(\alpha)$, containing the point $x$. It follows from \eqref{1_2016_04_14} that  
$$
\alpha^n_N(x)\supset\bigcap_{k=N}^n  f^{-k}\(B(f^k(x),e^{-k\beta})\). 
$$
Note also that 
$$
f^{-k}\(B(f^k(x),e^{-k\beta})\)\supset B\(x, e^{-k(\beta+\Delta)}\),
$$
where $e^\Delta$ is a Lipschitz constant, with respect to the spherical metric, of the map $f:\oc\to\oc$. Thus, 
$$
\alpha^n_N(x)
\supset \bigcap_{k=N}^n B\(x, e^{-k(\beta+\Delta)}\)
=B\(x, e^{-n(\beta+\Delta)}\)
$$
Finally, since 
$\alpha^n_0(x)=\alpha_0^{N-1}(x)\vee \alpha_N^n(x)$ and since $x\notin\bigcup_{n=0}^\infty\bigcup_{i\in I} 
f^{-n}(\partial U_i)$, there exists $\rho(x)>0$ such that $B(x,\rho(x))\subset \alpha_0^{N-1}(x)$. Thus, $\alpha^n_0(x)$ contains the ball $B\(x, e^{-n(\beta+\Delta)}\)$ for every integer $n\ge 1$ large enough (depending on $x$), and therefore $\alpha^n_0(x)\supset B\(x, C(x)e^{-n(\beta+\Delta)}\)$ for every integer $n\ge 1$, where $C(x)$ is some positive finite constant depending on $x$. Hence, the Weak Partition Existence Condition holds and property (iii) of Definition~\ref{dfn:Weakly_Markov} is established.
\end{proof}

\subsection{Conformal Graph Directed Markov Systems and Conformal IFSs}

\medskip 
Entering this subsection we start to deal with conformal systems. The ultimate difference between the examples to follow and those considered in the previous sections is that now we will be able to establish the convergence to the exponential one law, i.e. formulas \eqref{eq:wykl_LM1}--\eqref{eq:wykl_LM2B} for Full classes of radii and not merely $\beta$--Ultra Thick ones. Up to our best knowledge, this is the first time that the convergence to the exponential law is proved to hold for so general systems and measures along all radii.

In this subsection we apply our results about the exponential distribution of statistics of return times, namely Theorem~\ref{Main_Thm_2} for Weakly Markov systems and also the thin annuli property (Theorem~\ref{t:thin}, the same as Theorem~\ref{t:thin-Intro} from the introduction) for conformal IFSs to obtain Theorem~\ref{t1alr4} (the same as Theorem~\ref{t1alr4-Intro} from the introduction), i.e. the statistics of exponential one law for dynamical systems naturally induced by conformal GDMSs, in particular by conformal IFSs. So, let  
$$
\cS:=\{\phi_e:X_{t(e)}\to X_{i(e)}:e\in E\}
$$
be a conformal GDMS as defined in Section~\ref{CGDMS_Preliminaries} and let $A:E\times E\to\{0,1\}$ denote its incidence matrix. We assume throughout the subsection that $A$ (and so also $\cS$) is finitely irreducible. This time we however assume in addition that the Open Set Condition, in fact the Strong Open Set Condition of \cite{GDMS} holds. The Open Set Condition means that
\begin{equation}\label{osc}
\phi_a\left(\Int(X_{t(a)})\right)\cap \phi_b\left(\Int(X_{t(b)})\right)=\es
\end{equation}
whenever $a, b\in E$ with $a\ne b$. By a standard induction this condition implies that
\begin{equation}\label{osc_full}
\phi_\om\left(\Int(X_{t(\om)})\right)\cap \phi_\tau\left(\Int(X_{t(\tau)})\right)=\es
\end{equation}
whenever $\om$ and $\tau$ are any two incomparable words in $E_A^*$. The Strong Open Set Condition requires that in addition 
$$
J_\cS\cap \Int(X)\ne\es.
$$

\noindent Now let $f:E_A^\N\to\R$ be a H\"older continuous function, called in the sequel potential. We assume that $f$ is summable, meaning that
$$
\sum_{e\in E}\exp\(\sup(f|_{[e]})\)<+\infty.
$$
It is well known (see \cite{GDMS} or \cite{MU_Israel}) that the following limit 
$$
\P(f):=\lim_{n\to\infty}\frac1n\log\sum_{\om\in E_A^n}\exp\(\sup(f|_{[\om]})\)
$$
exists. It is called the topological pressure of $f$. It was proved in \cite{MU_Israel} (cf. \cite{GDMS}) that there exists a unique shift-invariant Gibbs/equilibrium measure $\mu_f$ for the potential $f$. The Gibbs property means that 
$$
C_f^{-1}\le\frac{\mu_f([\om|_n])}{\exp\(S_nf(\om)-\P(f)n\)}\le C_f
$$
with some constant $C_f\ge 1$ for every $\om\in E_A^\N$ and every integer $n\ge 1$, where here and in the sequel throughout this subsection
$$
S_n(g)=g_n(\om):=\sum_{j=0}^{n-1}g\circ \sg^j
$$
for every function $g:E_A^\N\to\C$. Let us record the following basic properties of the Gibbs state $\mu_f$.

\begin{fact}\label{Gibbs_Properties}
If the matrix $A$ is finitely irreducible and if $f:E_A^\N\to\R$ is a summable H\"older continuous potential, then the unique Gibbs state $\mu_f$ is ergodic and its topological support is equal to $E_A^\N$. In addition $\mu_f$ enjoys the Weak Independence Property {\rm (A)}.
\end{fact}

\fr Ergodicity has been proved in \cite{GDMS} while the Weak Independence Property (A) follows immediately from the definition of $\mu_f$. 

\medskip Following \cite{Urb2} we introduce the set
$$
\mathring{J}_\cS:=J_\cS\sms\bu_{\om\in E_A^*}\phi_\om(\partial X_{t(\om)}).
$$
We define 
$$
\mathring{E}_A^\N:=\pi_\cS^{-1}\(\mathring{J}\)
$$
and notice that for every $z\in\mathring{J}_\cS$ there exists a unique $\om(z)\in E_A^\N$ such that
$$
z=\pi(\om(z)).
$$
Moreover, $\om(z)\in \mathring{E}_A^\N$ and we simply denote it by $\pi^{-1}(z)$. Note that 
$$
\sg\(\mathring{E}_A^\N\)\sbt \mathring{E}_A^\N
$$
and this restricted shift map induces a map $T_\cS:\mathring{J}_\cS\to\mathring{J}_\cS$ by the formula
$$
T_\cS(z)=\pi\circ\sg(\pi^{-1}(z))\in\mathring{J}_\cS,
$$
so that the diagram 
\[\begin{CD}
\mathring{E}_A^\N @>\sigma>> \mathring{E}_A^\N\\
@V\pi VV @VV\pi V\\
\mathring{J}_\cS @>>T_\cS> \mathring{J}_\cS
\end{CD}
\]

\sp\fr commutes and the map $\pi:\mathring{E}_A^\N\to\mathring{J}_\cS$ is a continuous bijection. The map $T_\cS:\mathring{J}_\cS\to\mathring{J}_\cS$ is the main object of our interest in this subsection. Following notation of Section~\ref{sec:IFS} we denote
$$
\hat\mu_f:=\mu_f\circ\pi_\cS^{-1}.
$$
The following observation we deduce directly from Fact \ref{Gibbs_Properties}.

\begin{obs}\label{o2_2016_0211}
Suppose that $\cS$ is a finitely irreducible conformal GDMS satisfying the Strong Open Set Condition. If $f:E_A^\N\to\R$ is a summable H\"older continuous potential, then 
$$
\mu_f\(\mathring{E}_A^\N\)=1 
\  \and  \
\hat\mu_f\(\mathring{J}_\cS\)=1.
$$
Moreover, the projection $\pi:\mathring{E}_A^\N\to\mathring{J}_\cS$ establishes a measure--preserving isomorphism between measure--preserving dynamical systems $\(\sg:\mathring{E}_A^\N\to\mathring{E}_A^\N, \mu_f\)$ and $\(T_\cS:\mathring{J}_\cS\to\mathring{J}_\cS,\hat\mu_f\)$. 
\end{obs}

\fr We shall prove the following.

\begin{thm}\label{t1alr4}
Suppose that $\cS$ is a finitely irreducible and geometrically irreducible conformal GDMS satisfying the Strong Open Set Condition. If  $f:E_A^\N\to\R$ is a summable H\"older continuous potential such that
\begin{equation}\label{1_2016_02_22}
\sum_{e\in E}\exp\(\inf\(f|_{[e]}\)\)||\phi_e'||_\infty^{-\beta}<+\infty
\end{equation}
for some $\beta>0$, then the measure--preserving dynamical system $\(T_\cS:\mathring{J}_\cS\to\mathring{J}_\cS,\hat\mu_f\)$ is Weakly Markov and satisfies the Full Thin Annuli Property. In consequence, the exponential one laws of \eqref{eq:wykl_LM1} and \eqref{eq:wykl_LM2} hold for the dynamical system $\(T_\cS:\mathring{J}_\cS\to\mathring{J}_\cS,\hat\mu_f\)$. Precisely, 
\begin{equation}\label{eq:wykl_LM1C_3}
\lim_{r\to 0}	\sup_{t\geq0}\left|\hat\mu_f\left(\left\{z\in X\colon \tau_{B_r(x)}(z)>\frac{t}{\hat\mu_f(B_r(x))}\right\}\right)-e^{-t}\right|= 0
\end{equation}
for $\mu$--a.e. $x\in X$, i.e. the distributions of the normalized first entry time converge to the exponential one law, and 
\begin{equation}\label{eq:wykl_LM2C_3}
\lim_{r\to 0}	\sup_{t\geq0}\left|\hat\mu_{f{B_r(x)}}\left(\left\{z\in B_r(x)\colon \tau_{B_r(x)}(z)>\frac{t}{\hat\mu_f(B_r(x))}\right\}\right)-e^{-t}\right|= 0
\end{equation}
for $\mu$--a.e. $x\in X$, i.e.
the distributions of the normalized first return time converge to the exponential one law. 
\end{thm} 

\begin{proof}
Property (i) of being Weakly Markov (i.e. of Definition~\ref{dfn:Weakly_Markov}) for the dynamical sys\-tem $\(\sg:\mathring{E}_A^\N\to\mathring{E}_A^\N,\mu_f\)$ has been proved in \cite{GDMS}. For the dynamical system $\(T_\cS:\mathring{J}_\cS \to\mathring{J}_\cS,\hat\mu_f\)$ it then follows from the fact that the projection $\pi_\cS:E_A^\N\to J_\cS$ is H\"older continuous. Property (ii) has been also proved in \cite{GDMS}.  By virtue of Remark~\ref{rem:part}, in order to prove property (iii), it is enough to check that   the Weak Partition Existence Condition holds.  We do it now. Our proof resembles the one of property (iii) in Theorem~\ref{thm:ER}. We provide it for the sake of completeness. Let
$$
\alpha:=\{[e]\}_{e\in E}
$$
be the partition of $\mathring{E}_A^\N$ into cylinders of length one and let 
$$
\pi(\alpha):=\{\pi([e])\}_{e\in E}=\{\phi_e(J_\cS)\}_{e\in E}.
$$
Then
$$
\alpha_\sg^n=\{[\om]:\om\in E_A^n\}
$$
and
$$
\pi(\alpha)_T^n=\pi(\alpha_\sg^n)=\{\phi_\om(\mathring J_\cS):\om\in E_A^n\}.
$$
We know from \cite{GDMS} that $\h_{\mu_f}(\sg,\alpha)=\h_{\mu_f}(\sg)\in(0,+\infty)$, and so, by isomorphism, $\h_{\hat\mu_f}(T,\pi(\alpha))\in(0,+\infty)$. We also know from \cite{GDMS} that for 
$\mu_f$--a.e. $\om\in \mathring{E}_A^\N$, say $\om\in F\sbt\mathring{E}_A^\N$ with $\mu_f(F)=1$, the limit 
$$
\chi_{\mu_f}(\om)
:=-\lim_{n\to\infty}\frac1n\log\big|\phi_{\om|_n}'(\pi(\sg^n(\om)))\big|
$$
exists, is equal to 
$$
\chi_{\mu_f}
:=\int_{\mathring{E}_A^\N}\log\big|\phi_1'(\pi(\sg(\om)))\big|\,d\mu_f
$$
and belongs to $(0,+\infty)$. Fix $u\in V$, then fix $\xi\in \Int(X_u)$, and finally fix $R>0$ so small that 
$$
B(\xi,2R)\sbt \Int(X_u)
$$
Since $\mu_f\circ\pi_\cS^{-1}\(B(\xi,R)\)>0$, it follows from ergodicity of $\mu_f$ with respect to the shift map $\sg:E_A^{\N}\to E_A^{\N}$ and from Birkhoff's Ergodic Theorem that for $\mu_\psi$--a.e. $\om\in F$ there exists an infinite increasing sequence $(n_j)_{j=1}^\infty$ of positive integers such that
$$
\sg^{n_j}(\om)\in \pi_\cS^{-1}(B(\xi,R))
$$
for all $j\ge 1$ and 
$$
\lim_{j\to\infty}\frac{n_{j+1}}{n_j}=1.
$$
So, there exists a constant $A\ge 1$ such that
$$
n_{j+1}\le An_j
$$
for all $j\ge 1$. Now fix an arbitrary integer $k>n_1$. Then there exists a unique integer $j\ge 2$ such that
$$
n_{j-1}<k\le n_j.
$$
We then have $k>A^{-1}n_j$ Using the Distortion Property we conclude that, 
\begin{equation}\label{1_2016_09_10}
\begin{aligned}
\pi(\alpha)_T^k(\pi(\om))
&\spt\pi(\alpha)_T^{n_j}(\pi(\om))
=\phi_{\om|_{n_j}}\(\Int\(X_{t(\sg^{n_j}(\om))}\)\)
\spt \phi_{\om|_{n_j}}\(B(\xi,2R)\) \\
&\spt \phi_{\om|_n}\(B(\pi(\sg^{n_j}(\om)),R)\) \\
&\spt B\(\pi(\om),K^{-1}R\big|\phi_{\om|_{n_j}}'(\pi(\sg^{n_j}(\om)))\big|\) \\
&\spt B\(\pi(\om),\exp(-2\chi_{\mu_f}n_j)\)\\
&\spt B\(\pi(\om),\exp(-2A\chi_{\mu_f}k)\),
\end{aligned}
\end{equation}
where the last inclusion holds for all $n\ge 1$ large (depending on $\om$) enough, say $k\ge k(\om)\ge n_1$. On the other hand, obviously there exists some $\chi^*(\om)>0$ so large that
$$
\pi(\alpha)_T^k(\pi(\om)) \supset B\(\pi(\om),\exp(-\chi^*k)\)
$$
for all $k=0,1,\ldots, k(z)-1$. In conjunction with \eqref{1_2016_09_10} this gives that
$$
\pi(\alpha)_T^k(\pi(\om))
\supset B\Big(z,\exp\(-\max\big\{2A\chi_{\mu_f},\chi^*(z)\big\}k\)\Big)
$$
for all $k\ge 0$, and formula \eqref{borel_cantelli} is proved.    Since $\hat\mu_f(\pi(F))=1$, this establishes the Weak Partition Existence Condition. Property (A) trivially holds for $\hat\mu_f$ and (B) is satisfied because of \eqref{1_2016_02_22}. Since, see Theorem~\ref{t:thin}, measure $\hat\mu_f$ satisfies the Full Thin Annuli Property for IFSs, we are done in the case when $\cS$ is an IFS. In the general case we need an inducing argument. We only need to show that $\hat\mu_f$ satisfies the Full Thin Annuli Property. Fix $a\in E$ arbitrary and consider the following collection of $A$-admissible words.
$$
E_a:=\big\{\tau\in E_A^*:\tau_1=a,\, \forall_{2\le k\le|\tau|}\,\tau_k\ne a,\, A_{\tau_{|\tau|}a}=1\big\}.
$$
This gives rise to the following system of conformal uniformly contracting maps
$$
\cS_a:=\big\{\phi_\tau|_{\phi_a(X_{t(a)})}:\phi_a\(X_{t(a)}\)\to\phi_a\(X_{t(a)}\)\big\}.
$$
It is evident that $E_a^*\sbt E_A^*$ and that $\cS_a$ forms a conformal IFS whose limit set is contained in $J_\cS\cap \phi_a\(X_{t(a)}\)$; in the same vein the first return map $\sg_a:[a]\to[a]$ is canonically isomorphic to the full shift from $E_a^\N$ to $E_a^\N$. Moreover, the conditional measure $\mu_{a,f}$ on $[a]$ is the only Gibbs/equilibrium state of the shift map $\sg_a:[a]\to[a]$ and the H\"older continuous summable potential 
$$
E_a\ni\om\longmapsto S_{|\om_1|}f(\om)-\P(f)|\om_1|\in\R,
$$
where $|\om_1|$ maintains its original meaning as the length of a word in $E_A^*$ and $S_{|\om_1|}f$ denotes a Birkhoff's sum with respect to the original shift map $\sg:E_A^\N\to E_A^\N$. It therefore follows from the already proven cases of IFSs that $\mu_{a,f}$ satisfies the Full Thin Annuli Property. Since, by Poincar\'e Recurrence Theorem, 
$$
\mu_f\(\phi_e(X_{t(e)})\cap J_{\cS_e}\)=\mu_f\(\phi_e(X_{t(e)})\)
$$ 
and since 
$$
B(z,r)
\sbt \Int\(\phi_{\pi^{-1}(z)}(X_{t(\pi^{-1}(z))})\)
\sbt \phi_{\pi^{-1}(z)}\(X_{t(\pi^{-1}(z))}\)
$$
for all radii $r>0$ small enough, we therefore conclude that $\mu_f$ itself has the Full Thin Annuli Property. The proof is complete.
\end{proof}

\begin{rem}\label{r1_2026_0223}
Note that if the system $\cS$ of Theorem~\ref{t1alr4} is finite, then the hypothesis \eqref{1_2016_02_22} is automatically satisfied and can be removed from its assumptions. 
\end{rem}

Now we will pass to deal with measures that are of more geometric flavor. 

\begin{dfn}\label{fins}
We say that a real number $s$ belongs to $\Ga_\cS$, if
\begin{equation}\label{finite_parameters}
\sum_{e\in E}||\phi'_e||_\infty^s<+\infty.
\end{equation}
\end{dfn}

\fr We define the function $\zeta:E_A^\N\to\R$
$$
\zeta(\om):=\log\big|\phi_{\om_1}'(\pi_\cS(\sg(\om)))\big|.
$$
For every $t\in\R$ we consider the potential 
$$
t\zeta:E_A^\infty\to\R.
$$
Furthermore, we set
$$
\P(t):=\P(t\zeta).
$$

\noindent Let us record the following immediate observation.

\begin{obs}\label{o1_2016_01_20}
A real number $s$ belongs to $\Ga_\cS$ if and only if the H\"older continuous potential $s\zeta:E_A^\N\to\R$ is summable.
\end{obs}

\noindent We recall from \cite{MU_LMS} and \cite{GDMS} the following definitions:
$$
\g_\cS:=\inf\Ga_\cS=\inf\left\{s\in\R: \sum_{e\in E}||\phi_e'||_\infty^s<+\infty\right\}.
$$
Note that if the alphabet $E$ is finite, then $\g_\cS=-\infty$ and if $E$ is infinite, then $\g_\cS\ge 0$.
The proof of the following statement can be found in \cite{GDMS}.

\begin{prop}\label{p1_2016_01_12}
If $\cS$ is a finitely irreducible conformal GDMS, then for every $s\ge 0$ we have that 
$$
\Ga_\cS=\{s\in\R: \P(s)<+\infty\}
$$
In particular,
$$
\g_\cS=\inf\left\{s\in\R: \P(s)<+\infty\right\}.
$$
\end{prop}
For every $t\in \Ga_{\cS}$ we abbreviate
$$
\mu_t:=\mu_{t\zeta}.
$$
As an immediate consequence of Theorem~\ref{t1alr4} we get the following.

\begin{cor}\label{c1alr7}
Suppose that $\cS$ is a finitely irreducible and geometrically irreducible conformal GDMS satisfying the Strong Open Set Condition. Fix a real number $t>\g_\cS$. Then the corresponding measure--preserving dynamical system $\(T_\cS:\mathring{J}_\cS\to\mathring{J}_\cS,\hat\mu_t\)$ is Weakly Markov and satisfies the Full Thin Annuli Property. In particular, the exponential one laws of \eqref{eq:wykl_LM1} and \eqref{eq:wykl_LM2} hold for the dynamical system $\(T_\cS:\mathring{J}_\cS\to\mathring{J}_\cS,\hat\mu_t\)$.
\end{cor}

\begin{rem}\label{r1_2016_02_23}
Recall that if the system $\cS$ of Corollary~\ref{c1alr7} is finite, then $\gamma_{\cS}=-\infty$ and the hypothesis $t>\gamma_\cS$ is automatically fulfilled.
\end{rem}

\begin{rem}\label{r2_2016_10_19}
In the setting of Corollary~\ref{c1alr7}, let $h_\cS$ be the Hausdorff dimension of the limit set $J_\cS$. It is known, see \cite{GDMS}, that then $\H_{h_\cS}(J_\cS)$, the $h_\cS$--dimensional Hausdorff measure of $J_\cS$ is finite while the corresponding packing measure $\P_{h_\cS}(J_\cS)$ is positive. If either one of these two measures is both finite and positive, then this measure is equivalent to the measure $\hat\mu_{h_\cS}$ (which then does exist!) with uniformly bounded Radon--Nikodym derivatives. Thus the Full Thin Annuli Property holds respectively for $\H_{h_\cS}$ or $\P_{h_\cS}$ (or both) restricted to $J_\cS$. This is always the case when the system $\cS$ is finite. Note also that if $\mu$ is any finite Borel measure satisfying Ahlfors property with exponent $h>d-1$, then, as a 
straight volume argument shows, this measure satisfies the Full Thin Annuli Property at each point of its topological support.
\end{rem}

Since, as it is well known, the harmonic measure of the limit set of a finite conformal IFS which satisfies the strong separation condition, is equivalent, with uniformly bounded Radon--Nikodym derivatives, to a Gibbs/equilibrium measure, as an immediate consequence of Theorem~\ref{t1alr4}, we get the following.

\begin{cor}\label{c1alr7new}
Suppose that $\cS$ is a conformal IFS in the complex plane $\C$ satisfying the Strong Separation Condition. Then the harmonic measure of its limit set satisfies the Full Thin Annuli Property.
\end{cor}

\fr In fact this corollary is a consequence of Theorem~\ref{t1alr4} under the additional assumption of geometrical irreducibility. In the real--analytic case, still IFS, it follows, from easy to prove, upper estimates of the harmonic measure of a ball by its radius raised to some positive power. 

\subsection{Conformal Parabolic GDMSs} 
In this subsection, following \cite{MU_Parabolic_1}
and \cite{GDMS}, we first shall provide the appropriate setting and basic properties of conformal parabolic iterated function systems, and more generally of parabolic graph directed Markov systems. We then prove for them the appropriate theorems on convergence to the exponential law. 

\sp\fr As in Section~\ref{CGDMS_Preliminaries} there are given a directed multigraph $(V,E,i,t)$ ($E$ and $V$ both (!) finite), an incidence matrix $A:E\times E\to \{0,1\}$, and two functions $i,t:E\to V$ such that $A_{ab} = 1$ implies $t(b) = i(a)$. Also, we have nonempty compact metric spaces $\{X_v\}_{v\in V}$ and their respective bounded connected neighborhoods $W_v$, $v\in V$. Suppose further that we have  a collection of conformal maps $\phi_e:X_{t(e)}\to X_{i(e)}$, $e\in E$, satisfying the following conditions:

\sp\begin{itemize}
\item[(1)]\emph{Open Set Condition:}
 $\phi_i(\Int(X))\cap \phi_j(\Int(X))=\es$ for all $i\ne j$.

\sp\item[(2)] $|\phi_i'(x)|<1$ everywhere except for finitely many
pairs $(i,x_i)$, $i\in E$, for which $x_i$ is the unique fixed point
of $\phi_i$ and $|\phi_i'(x_i)|
=1$. Such pairs and indices $i$ will be called parabolic and the set of
parabolic indices will be denoted by $\Om$. All other indices will be 
called hyperbolic. We assume that $A_{ii}=1$ for all $i\in\Om$.

\sp\item[(3)]  $\forall n\ge 1 \  \forall \om = (\om_1,...,\om_n)\in E_A^n$
if $\om_n$ is a hyperbolic
index or $\om_{n-1}\ne \om_n$, then $\phi_{\om}$ extends conformally to
an open connected set $W_{t(\om_n)}\sbt\R^d$ and maps $W_{t(\om_n)}$ into $W_{i(\om_n)}$.

\sp\item[(4)] If $i$ is a parabolic index, then $\bi_{n\ge
0}\phi_{i^n}(X)
=\{x_i\}$ and the diameters of the sets $\phi_{i^n}(X)$ converge
to 0.

\sp\item[(5)] $\exists s<1 \  \forall n\ge 1  \  \forall \om\in E_A^n$ if 
$\om_n$ is a hyperbolic index or $\om_{n-1}\ne \om_n$, then
$$
\|\phi_\om'\|\le s.
$$
\end{itemize}

\sp\fr We call such a system of maps 
$$
\cS=\{\phi_i:i\in E\}
$$ 
a subparabolic
iterated function system. Let us note that conditions (1), (3), (5) are modeled on similar  conditions which were used to examine hyperbolic
conformal systems. If $\Om\ne\es$, we call the system  
$\{\phi_i:i\in E\}$ parabolic. As
declared in (2) the elements of the set $E\sms \Om$ are called
hyperbolic.
We extend this name to all the words appearing in (5). It follows
from (3) that for every hyperbolic word $\om$,
$$
\phi_\om(W_{t(\om)})\sbt W_{t(\om)}.
$$
Note that our conditions ensure that $\phi_i'(x) \neq 0$ 
for all $i\in E$ and all $x \in X_{t(i)}$. It was proved (though only for IFSs but the case of GDMSs can be treated completely similarly) in 
\cite{MU_Parabolic_1} (comp. \cite{GDMS}) that
\begin{equation}\label{1_2016_03_15}
\lim_{n\to\infty}\sup_{|\om|=n}\{\diam(\phi_\om(X_{t(\om)}))\}=0.
\end{equation}
As its immediate consequence, we record the following.

\begin{cor}\label{p1c2.3} 
The map $\pi:E_A^\infty\to X:=\bigoplus_{v\in V}X_v$, $$
\{\pi(\om)\}:=
\bi_{n\ge 0}\phi_{\om|_n}(X),
$$
is well defined, i.e. this intersection is always a singleton, and the map $\pi$ is uniformly continuous.
\end{cor}

\fr As for hyperbolic (attracting) systems the limit set $J = J_\cS$ of the system $\cS = \{\phi_e\}_{e\in e}$ is defined to be
$$
J_\cS:=\pi(E_A^\infty)
$$
and it enjoys the following self-reproducing property:
$$
J = \bu_{e\in E} \phi_e(J).
$$
We now, following still \cite{MU_Parabolic_1} and \cite{GDMS}, want to associate to the parabolic system $\cS$ a canonical hyperbolic system $\cS^*$. The set of edges is this.
$$
E_*:= \big\{i^nj: n\ge 1, \  i\in\Om, \ i\ne j\in E, \ A_{ij}=1\big\} \cup 
(E\sms \Om)\sbt E_A^*.
$$ 
We set
$$
V_*=t(E_*)\cup i(E_*)
$$
and keep the functions $t$ and $i$ on $E_*$ as the restrictions of $t$ and $i$ from $E_A^*$. The incidence matrix $A_*:E_*\times E_*\to\{0,1\}$ is defined in the natural (the only reasonable) way by declaring that $A^*_{ab}=1$ if and only if $ab\in E_A^*$. Finally 
$$
\cS^*:=\{\phi_e:X_{t(e)}\to X_{t(e)}:\, e\in E_*\}.
$$
It immediately follows from our assumptions (see \cite{MU_Parabolic_1} and \cite{GDMS} for details) that the following is true.

\begin{thm}\label{p1t5.2} 
The system $\cS^*$ is a hyperbolic conformal GDMS and the limit sets $J_\cS$ and $J_{\cS^*}$ differ only by a countable set.
\end{thm}

\fr We have the following quantitative result, whose complete proof can be found in \cite{ADU}. 

\begin{prop}\label{p1c5.13} 
Let $\cS$ be a conformal parabolic GDMS. Then there exists a constant $C\in(0,+\infty)$ and for every $i\in\Om$ there exists some constant
$\beta_i\in(0,+\infty) $ such that for all $n\ge 1$ and for all $z\in X_i:=
\bu_{j\in I\sms\{i\}}\phi_j(X)$,
$$
C^{-1}n^{-\frac{\beta_i+1}{\beta_i}}\le |\phi_{i^n}'(z)|
\le Cn^{-\frac{\beta_i+1}{\beta_i}}.
$$
In fact we know more: if $d=2$ then all constants $\b_i$ are integers $\ge 1$ and if $d\ge 3$, then all constants $\b_i$ are equal to $1$. 
\end{prop}

\fr Let 
$$
\b=\b_\cS:=\min\{\b_i: i\in\Om\}.
$$
Passing to equilibrium/Gibbs states and their escape rates, we now describe the class of potentials we want to deal with. This class is somewhat narrow as we restrict ourselves to geometric potentials only. There is no obvious natural larger class of potentials for which our methods would work and trying to identify such classes would be of dubious value and unclear benefits. We thus only consider potentials of the form
$$
E_A^\infty\ni\om\mapsto\zeta_t(\om)
:=t\log\big|\phi_{\om_0}'(\pi_\cS(\sg(\om)))\big|\in\R, \  \  t\ge 0.
$$
We then define the potential $\zeta_t^*:E_{*A^*}^\infty\to\R$ as
$$
\zeta_t^*(i^nj\om)=\sum_{k=0}^n\zeta_t(\sg^k(i^nj\om)), \  \  \ i\in\Om, \  n\ge 0, \  j\ne i \  \and  \  i^nj\om\in  E_{*A^*}^\infty.
$$
We shall prove the following.

\begin{prop}\label{p1ps1}
If $\cS$ is a finite conformal parabolic GDMS, then given $t\ge 0$ the potential $\zeta_t^*$ is H\"older continuous. Moreover, this potential is summable if and only if 
$$
t>\frac{\b}{\b+1}.
$$
\end{prop}

\bpf
H\"older continuity of potentials $\zeta_t^*$, $t\ge 0$, follows from the fact that the system $\cS^*$ is hyperbolic, particularly from its distortion property, while the summability statement immediately follows from Proposition~\ref{p1c5.13}.
\epf

\sp\fr So, for every $t>\frac{\b}{\b+1}$ we can define $\mu_t^*$ to be the unique equilibrium/Gibbs state for the potential $\zeta_t^*$ with respect to the shift map $\sg_*:
E_{*A^*}^\infty\to E_{*A^*}^\infty$. We know that $\mu_t^*$ gives rise to a Borel $\sg$-finite, unique up to multiplicative constant, $\sg$-invariant measure $\mu_t$ on $E_A^\infty$, absolutely continuous, in fact equivalent, with respect to $\mu_t^*$; see \cite{GDMS} for details in the case of $t=b_\cS=b_{\cS^*}$, the Bowen's parameter of the systems $\cS$ and $\cS^*$ alike. The case of all other $t>\frac{\b}{\b+1}$ can be treated similarly. It follows from
\cite{GDMS} that the measure $\mu_t$ is finite if and only if either

\sp\begin{itemize}
\item[(a)] $t\in\left(\frac{\b}{\b+1},b_\cS\right)$ or

\sp\item[(b)] $t=b_\cS$ \and $b_\cS>\frac{2\b}{\b+1}$. 
\end{itemize}

\sp\fr The main result of this subsection is the following.

\begin{thm}\label{c1alr7B}
Suppose that $\cS$ is a finite irreducible and geometrically irreducible parabolic conformal GDMS satisfying the Strong Open Set Condition. Fix a real number $t$ for which one of the conditions (a) or (b) above holds. Then the corresponding measure--preserving dynamical system $\(T_\cS:\mathring{J}_\cS\to\mathring{J}_\cS,\hat\mu_t\)$ satisfies the Full Thin Annuli Property and the exponential one laws of \eqref{eq:wykl_LM1} and \eqref{eq:wykl_LM2} hold for the dynamical system $\(T_\cS:\mathring{J}_\cS\to\mathring{J}_\cS,\hat\mu_t\)$.
\end{thm}

\bpf
The proof consists of the following three ingredients. The first one is that the induced system $\cS^*$ with the measure $\mu_t^*$ satisfies all the hypotheses of Corollary~\ref{c1alr7}. The second one is that the measure--preserving dynamical system $\(T_{\cS^*}:\mathring{J}_{\cS^*}\to\mathring{J}_{\cS^*},\hat\mu_t^*\)$ forms the 1st return time map of the measure--preserving system $\(T_\cS:\mathring{J}_\cS\to \mathring{J}_\cS,\hat\mu_t\)$. The third one is that, according to one of the main results of \cite{HSV}, if the 1st return time map of a measure--preserving dynamical systems satisfies the exponential one laws of \eqref{eq:wykl_LM1} and \eqref{eq:wykl_LM2}, then so does the original system. 
\epf 

We would like to remark that Theorem~\ref{c1alr7B} covers such examples as Parabolic Cantor Sets (see \cite{Ur}) Apollonian packing system (see \cite{GDMS}), and finitely generated Schottky groups with some generating ball tangent to each other, Farey map, and much more. More information about these systems can be found for example in \cite{GDMS}.

\sp We would like however to single out one class of parabolic systems, namely parabolic rational functions. These are defined as rational functions of the Riemann sphere whose restrictions to their Julia sets are expansive but not expanding. Equivalently (see \cite{DU_2}), those whose Julia sets contain no critical points but do contain rationally indifferent (parabolic) periodic points. Such rational functions admit Markov partitions with arbitrarily small diameters (see \cite{DU_2} again). Thus, these can be viewed as finite parabolic conformal iterated function systems. Let $f:\hat\C\to\hat\C$ be such a function. Let $h_f$ be the Hausdorff dimension of the Julia set $J(f)$ of $f$. Let $p\ge 1$ denote the maximal number of petals around parabolic periodic points of $f$. It coincides with the number $\beta$ of the above mentioned parabolic iterated function system. Suppose that
\begin{equation}\label{1_2016_10_19}
h_f>\frac{2p}{p+1}.
\end{equation}
We know from \cite{DU_2}, \cite{DU_3}, and \cite{ADU} that if \eqref{1_2016_10_19} holds, in fact if $h_f\ge 1$, then the $h_f$-dimensional Hausdorff measure of $J(f)$ is positive and finite. Furthermore (see the same three papers), still assuming \eqref{1_2016_10_19},
there exists then a unique probability $f$-invariant measure absolutely continuous, in fact equivalent, with respect to this Hausdorff measure.  This invariant measure coincides with the measure $\hat\mu_{h_f}$ obtained from the parabolic conformal iterated function system generated by the above mentioned Markov partition. Therefore, Theorem~\ref{c1alr7B} entails the following.

\begin{thm}\label{c1alr7c}
Let $f:\hat\C\to\hat\C$ be a parabolic rational function whose Julia sets is not contained in any real analytic curve. Assume also that \eqref{1_2016_10_19} holds. Then the measure $\hat\mu_{h_f}$  satisfies the Full Thin Annuli Property and the exponential one laws of \eqref{eq:wykl_LM1} and \eqref{eq:wykl_LM2} hold for the dynamical system $(f:J(f)\to J(f),\hat\mu_{h_f}\)$.
\end{thm}

\begin{rem}
A classical example to which Theorem~\ref{c1alr7c} applies is the polynomial $\hat\C\ni z\mapsto z^2+\frac14$.
\end{rem}
 
\begin{rem}
The obvious analogue of Theorem~\ref{c1alr7c} holds for all $t\in \left(\frac{p}{p+1}, h_f\right)$. Note however that this case is also covered by Subsection~\ref{holder_equil}.
\end{rem}

\subsection{Conformal Expanding Repellers}\label{CER}

\sp Now let us formulate the definition of a conformal expanding repeller, the primary object of interest in this subsection.

\begin{dfn}\label{exprep}
Let $U$ be an open subset of $\R^d$, $d\ge 1$. Let $J$ be a compact subset of $U$. Let $T:U\to\R^d$ be a conformal map. The map $T$ is called a conformal expanding repeller if the following conditions are satisfied:
\begin{enumerate}
\item{}$T(J)=J$,

\sp\item{} 
$|T'_{|J}|>1$,

\sp\item{} there exists an open set  $V$ such that $\overline{V}\subset U$ and
$$
J=\bigcap_{k=0}^\infty T^{-n}(V).
$$
\item{}the map $T_{|J}:J\to J$ is topologically transitive.
\end{enumerate}
\end{dfn}

\fr So, a conformal expanding repeller is an expanding repeller of Subsection~\ref{DEMER} for which the corresponding map $T$ is not merely smooth but conformal. Typical examples of conformal expanding repellers are provided by the following. 

\begin{prop}\label{p1_2015_05_11} 
If $f:\hat\C\to\hat\C$ is a rational function of degree $d\ge 2$, such that the map $f$ restricted to its Julia set $J(f)$ is expanding, then $J(f)$ is a conformal expanding repeller.
\end{prop}

\begin{thm}\label{t1_2016_10_20}
Let $T:J\to J$ be a conformal expanding repeller such that $J$ is not contained in any real analytic submanifold of dimension $\le d-1$. Let $\psi:J\to\R$ be a H\"older continuous potential and, see \cite{PUbook}, let $\mu_\psi$ be the corresponding equilibrium (also frequently referred to as Gibbs) state. Then the measure--preserving dynamical system $\(T:J\to J,\mu_\psi\)$ is Weakly Markov and satisfies the Full Thin Annuli Property. In particular, the exponential one laws of \eqref{eq:wykl_LM1} and \eqref{eq:wykl_LM2} hold for the dynamical system $\(T:J\to J,\mu_\psi\)$. 
\end{thm}

\begin{proof}
The dynamical system $\(T:J\to J,\mu_\psi\)$ is Weakly Markov because this property was established in Theorem~\ref{thm:ER} for all expanding repellers. For the Full Thin Annuli Property we use again Markov partitions discussed and utilized in Subsection~\ref{DEMER}, particularly in the proof of Theorem~\ref{thm:ER}. So, let 
$$
\cR=\{R_1,\ldots,R_q\}
$$
be the Markov partition considered therein. We now associate to $\cR$ a finite conformal graph directed Markov system. The set of vertices is equal to $\cR$ while the alphabet $E$ is defined as follows.
$$
E:=\big\{(i,j)\in \{1,2\,\ldots,q\}\times \{1,2\,\ldots,q\}:\Int R_j\cap T(\Int R_i) \ne \es \big\}. 
$$
Now, for every $(i,j)\in E$ there exists a unique conformal map $T_{i,j}^{-1}:B(R_j,\delta)\to\R^d$ such that
$$
T_{i,j}^{-1}(R_j)\sbt R_i.
$$
Define the incidence matrix $A:E\times E\to\{0,1\}$ by
$$
A_{(i,j) (k,l)}=
\begin{cases}
1 \ &\text{{\rm  if }} \ l=i \\
0 \ &\text{{\rm  if }} \ l\ne i.
\end{cases}
$$
Define further
$$
t(i,j)=j \  \ \text{{\rm and }} \  \  i(i,j)=i. 
$$
Of course 
$$
\cS_\cR=\{T_{i,j}^{-1}:(i,j)\in E\}
$$ 
forms a finite conformal directed Markov system, and $\cS_\cR$ is irreducible since the map $T:J\to J$ is transitive. In addition, $\cS_\cR$ is geometrically irreducible because $J$ is not contained in any any real analytic submanifold of dimension $\le d-1$. Define the potential $\hat\psi:E_A^\N\to\R$ by 
$$
\hat\psi:=\psi\(\pi_{\cS_\cR}\).
$$
The potential $\hat\psi$ is H\"older continuous as a composition of two H\"older continuous functions; H\"older continuity of $\pi_{\cS_\cR}$ with a standard metric on the symbol space follows immediately from the expanding property and a detailed proof can be found e.g. in \cite{PUbook}. Moreover, it is known (see e.g. \cite{PUbook}) that
$$
\mu_\psi=\mu_{\hat\psi}\circ\pi_{\cS_\cR}^{-1}.
$$
Therefore, the Full Thin Annuli Property of $\mu_\psi$ follows from Theorem~\ref{t1alr4}; remember that this is not a property of a system but of a measure. The proof is complete.
\end{proof}

\begin{rem}\label{r2_2016_10_19B}
Since, every conformal expanding repeller $T:J\to J$ admits a finite Markov partition (see Theorem~\ref{t1_2016_10_20}), the proof of this theorem shows that Corollary~\ref{c1alr7} and Remark~\ref{r2_2016_10_19} now apply, the latter for IFSs with finite alphabets. These two thus yield the Full Thin Annuli Property of $\HD(J)$--dimensional Hausdorff measure of $J$.
\end{rem}

\subsection{Equilibrium States for Rational Maps of the Riemann Sphere $\oc$ and H\"older Continuous Potentials with a Pressure Gap} \label{holder_equil}

Let $f:\oc\to\oc$ be a rational map of degree larger than $1$. Denote by  $J(f)$ its Julia set. Let $\varphi:J(f)\to\mathbb{R}$ be a H\"older
continuous function. As in previous subsections keep $\P(\varphi)$ to denote its topological pressure with respect to the dynamical system generated by the map $f:J(f)\to J(f)$. M. Lyubich proved in \cite{lyubich} that in our context of rational functions each continuous function admits an equilibrium  state. It was shown in \cite{dueqrat} that if $\varphi$ (being H\"older
continuous) has a pressure gap, i.e. if
$$
\P(\varphi)>\frac1n\sup(S_n\varphi)
$$
for some integer $n\ge 1$, then there exists a unique equilibrium measure for $\varphi$ which we again denote by $\mu_\varphi$.

In \cite{suz} several strong stochastic properties of this equilibrium measure $\mu_\varphi$ have been deduced from a special inducing scheme. The induced map forms a conformal Iterated Function System, satisfying the Strong Separation Condition, in particular the Strong Open Set Condition.
 Now we shall prove the following main result of this section.

\begin{thm}\label{thm: equil}
Let $f:\oc\to\oc$ be an arbitrary rational map of degree larger
than $1$ whose Julia set is not contained in a real analytic curve (this is always the case if for instance $\HD(J(f))>1$). Let $\varphi:J(f)\to\mathbb{R}$ be a H\"older
continuous function with a pressure gap. Then $(J(f), f, \mu_\varphi)$ forms a Weakly Markov system with the Full Thin Annuli Property.
Consequently, the exponential one laws of \eqref{eq:wykl_LM1} and \eqref{eq:wykl_LM2} hold for the dynamical system $\(T:J(f)\to J(f),\mu_\varphi\)$. 
\end{thm}

\begin{proof}
In order to  check that the required properties hold, we refer to appropriate  results in \cite{suz}. The argument for the dynamical system $\(T:J(f)\to J(f),\mu_\varphi\)$ to be Weakly Markov is actually the same as the one presented in the proof of Theorem~\ref{thm: projective_equil}. It is exactly the same for items (i) and (iii) while the argument for 
item (ii) is simpler; it holds in the current one-dimensional setting, since the limit under consideration exists $\mu_\varphi$--a.e. and is equal to the Hausdorff dimension of the measure $\mu_\varphi$, which is a positive number.

\sp The Full Thin Annuli Property is a consequence of the above mentioned fine inducing procedure, see 
\cite{suz}, Section ~3. We follow the notation of \cite{suz}, especially Section~8 of this paper. The fine inducing construction leads to a conformal Iterated Function System, satisfying the Strong Separation Condition, and such that the limit set of this system is of full $\mu_\varphi$ measure. We denote this system by $\cS$.
We recall briefly the way this induced system is constructed. For a properly chosen topological disc $U$, the system $\cS$ is defined by a family of conformal univalent homeomorphisms $\phi_e:U\to D_e$, $e\in E$, where $E$ is some countable set and  $\overline D_e\subset U$ for every $e\in E$. Each map  $\phi_e$, $e\in E$, is  just, a suitably chosen, holomorphic branch of the inverse of some iterate of $f$, say $f^{N(e)}$, mapping $U$ onto $D_e$. As usual, denote the corresponding projection from $E^\N$ to $\oc$ by $\pi_\cS$. The iterated function system $\cS$, together with the summable H\"older potential 
$$
\overline\varphi=S_{N(e)}\varphi\circ\pi_\cS-\P(\varphi)N(e):E^\N\to\R,
$$
arising naturally from the inducing procedure, admits an (invariant) equilibrium state which is equivalent to the initial measure 
$\mu_\varphi$. We claim that the IFS $\cS$ together with the (induced) potential $\overline\varphi$, satisfies the hypotheses of Theorem~\ref{t1alr4}, with $f$  therein being replaced by $\overline\varphi$. 
We shall  sketch the argument here, referring to appropriate estimates in \cite{suz}.
The estimate which we need to verify the assumption of Theorem~\ref{t1alr4} is the following (see \eqref{1_2016_02_22})
\begin{equation}\label{eq:summ}
\sum_{e\in E}\exp\(\inf\(\overline\varphi|_{[e]}\)\)||\phi_e'||_\infty^{-\beta}<+\infty
\end{equation}
with some $\beta>0$. Note, however, that $\exp\(\inf\(\overline\varphi|_{[e]}\)\)$ is multiplicatively comparable to $\mu_\varphi(D_e)$ independently of $e$, and, consequently, in order to verify \eqref{eq:summ}, it is enough to check that 
\begin{equation}\label{eq:integrable}
\int |F'|^\beta d\mu_\varphi<\infty \quad  \text {for some}\quad  \beta>0,
\end{equation}
where the map $F$ is defined on each set $D_e$ just as $(\phi_e)^{-1}$.
This can be easily done by using the estimates provided in \cite{suz}. Indeed, by the definition of $F$ and the system $\cS$, we have that $F|_{D_e}=f^{N(e)}|_{D_e}$.  Moreover, the estimates in \cite{suz} (see e.g. the formula (3.1) in 
\cite{suz}) show that 
$$
\mu_{\varphi}\bigg(\bigcup_{e\colon N(e)\ge n} D_e\bigg)\le 2e^{-n\gamma}
$$ 
for every integer $n\ge 1$ and some $\gamma>0$. Using the trivial estimate $|F'|\le ||f'||^{N(e)}$, \eqref{eq:integrable} follows immediately. Finally, the system $\cS$ is geometrically irreducible since the Julia set is not contained in a real analytic curve. 

Therefore, we are in position to apply Theorem~\ref{t1alr4}, and the  measure $\mu_\varphi$ has the Full Thin Annuli Property. The proof is complete.
\end{proof}

\subsection{Dynamically Semi--Regular Meromorphic Functions}
Let $f:\C\to \oc$ be a    meromorphic function. Let $\Sing(f^{-1})$ be the set of all singular points of $f^{-1}$, i.e. the set of all points $w\in\oc$ such that if $W$ is any open connected neighborhood of $w$, then there exists a connected component $U$ of $f^{-1}(W)$ such that the map $f:U\to W$ is not bijective. Of course, if $f$ is a rational function, then $\Sing(f^{-1})=f(\Crit(f))$. Define
$$
\PS(f):=\bu_{n=0}^\infty f^n(\Sing(f^{-1})).
$$
The function $f$ is called \emph{topologically hyperbolic} if
$$
\dist_{\text{Euclid}}(J_f ,\PS(f)) >0,
$$
and it is called \emph{expanding} if there exist $c>0$ and $ \lambda>1$
such that
$$
|(f^n)'(z)|\ge c\lambda^n 
$$
for all integers $n\ge 1$ and all points $z\in J_f\sms
f^{-n}(\infty)$. Note that every topologically hyperbolic meromorphic
function is \emph{tame} (see definition before Theorem \ref{prop:1}).  A meromorphic function that is both topologically
hyperbolic and expanding 
is called \emph{hyperbolic}. The meromorphic function $f:\C\to\oc$ is
called dynamically {\it semi-regular} if it is of finite order, commonly
denoted by $\rho_f$, and satisfies the following rapid
growth condition for its derivative.  
\begin{equation}\label{eq intro}
|f'(z)|\geq \kappa ^{-1} (1+|z|)^{\alpha _1} (1+|f(z)|)^{\alpha_2} \; ,
\quad z\in J_f, 
\end{equation}
with some constant $\kappa >0$ and $\alpha_1,\alpha_2$ such that $\alpha_2 > \max\{-\alpha_1 ,0\}$.  Set $\alpha:=\alpha_1+\alpha_2$.

\ 

\fr Let $h:J_f\to\R$ be a weakly
H\"older continuous function 
in the sense of \cite{MayUrb10}. The definition, introduced in
\cite{MayUrb10} is somewhat technical and we will not provide it in
the current paper. What is important is that each bounded, uniformly
locally H\"older function $h:J_f\to\R$ is weakly H\"older. Fix
$\tau>\alpha_2$ as required in \cite{MayUrb10}. For $t\in\R$, let
\begin{equation}
  \label{eq:7}
  \psi_t=-t\log|f'|_\tau+h
\end{equation}
where $|f'(z)|_\tau$ is the norm, or, equivalently, the scaling
factor, of the derivative of $f$ evaluated at a point $z\in J_f$ with
respect to the Riemannian metric
$$
|d\tau(z)|=(1+|z|)^{-\tau}|dz|.
$$
For any $t>\rho_f/\alpha$ let $\cL_t:C_b(J_f)\to C_b(J_f)$ be the
corresponding \emph{Perron--Frobenius operator} given by the formula
$$
\cL_tg(z)=\sum_{w\in f^{-1}(z)}g(w)e^{\psi_t(w)}.
$$
The hypothesis $t>\rho_f/\alpha$ guaranties that the series
$$
\sum_{w\in f^{-1}(z)}|f'(w)|_\tau^{-t}
$$
converges uniformly on $J_f$, and, in particular, the linear operator $\cL_t:C_b(J_f)\to C_b(J_f)$ is well defined and bounded.
It was shown in \cite{MayUrb10} that, for every $z\in J_f$
, the limit
$$
\lim_{n\to\infty}\frac1n\log\cL_t\1(z)
$$
exists and takes on the same common value, which we denote by $\P(t)$
and call \emph{the topological pressure} of the potential
$\psi_t$. The following theorem was proved in \cite{MayUrb10}.

\

\begin{thm}\label{t1dns111}
If $f:\C\to\oc$ is a dynamically semi-regular meromorphic function and
$h:J_f\to\R$ is a weakly H\"older continuous potential, then for every
$t>\rho_f/\alpha$ there exist uniquely determined Borel probability measures $m_t$ and $\mu_t$ (which do depend on the function $h$ too even though this is not explicitly indicated) on $J_f$ with the following properties.

\sp\begin{itemize}
\item[{\rm(a)}] \ $\cL_t^*m_t=m_t$.

\sp\item[{\rm(b)}] \ $\P(t)=\sup\big\{\h_\mu(f)+\int\psi_t\, d\mu:\mu\circ f^{-1}=\mu \  \
  \text{{\rm and }}  \ \int\psi_t\, d\mu>-\infty\big\}$.

\sp\item[{\rm(c)}] \ $\mu_t\circ f^{-1}=\mu_t$, $\int\psi_t\, d\mu_t>-\infty$, \  and  \
$
\h_{\mu_t}(f)+\int\psi_t\,d\mu_t=\P(t).
$

\sp\item[{\rm(d)}] \ The measures $\mu_t$ and $m_t$ are equivalent and the
  Radon--Nikodym derivative $\frac{d\mu_t}{dm_t}$ has a
  nowhere-vanishing H\"older continuous version which is bounded from above.
\end{itemize}
\end{thm}

\fr Item (a) (along with (d)) essentially means that $m_t$ and $\mu_t$ are Gibbs states of the potential $\psi_t$, while items (b) and (c) mean that $\mu_t$ is an equilibrium state for the potential $\psi_t$. We shall prove the following.

\begin{thm}\label{t1_2016_02_22}
Let $f:\C\to\oc$ be a dynamically semi-regular meromorphic function whose Julia set is not contained in a real analytic curve (this is always the case if for instance $\HD(J_f)>1$). Let $t>\rho_f/\alpha$, and let $h:J_f\to\R$ be a weakly H\"older continuous potential. Then the measure--preserving dynamical system $\(f:\C\to\oc,\mu_t\)$ is Weakly Markov and satisfies the thin annuli property. In particular, the exponential one laws of \eqref{eq:wykl_LM1} and \eqref{eq:wykl_LM2} hold for the dynamical system $\(f:\C\to\oc,\mu_t\)$.
\end{thm} 

\begin{proof}
Property (i) of being Weakly Markov (i.e. of Definition~\ref{dfn:Weakly_Markov}) for the dynamical system $\(f:\C\to\oc,\mu_t\)$ has been proved in \cite{MayUrb10} as Theorem~6.16. Property (ii) is a part of Theorem~8.1 therein. For property (iii) it suffices to notice that the Weak Partition Existence Condition holds. And it does because of the first displayed formula after (6.22) in the proof of Theorem~6.25 (Variational Principle) in \cite{MayUrb10}, and because the map $f$ is expanding. 

\sp We are thus left to prove the Thin Annuli Property. As in the case of conformal graph directed Markov systems it will be based on an inducing argument. The point is that one can construct conformal IFSs having any given non-periodic recurrent point of the Julia set in the interior of its seed set. We formulate the appropriate theorem in a more general setting which does not enlarge the volume of our considerations. Following \cite{PU_tame} and \cite{SkU} we call a meromorphic function $f:\C\to\oc$ \emph{tame} if
$$
J(f)\sms \overline{\PS(f)}\ne\es.
$$
The following theorem was proved in \cite{Dob11}. 
 
\begin{thm}\label{prop:1}
Let $f:\C\to\oc$ be a tame meromorphic function. Fix a non-periodic point $z\in J(f)\sms
\overline{\PS(f)}$, $\kappa>1$, and $K>1$. Then for all $\lambda>1$ and 
for all $r>0$ sufficiently small there exists
an open connected set $V=V(z,r)\sbt\C\sms\overline{\PS(f)}$ such that
  \begin{itemize}
  \item[(a)] If $U\in \Comp(f^{-n}(V))$ and $U\cap V\neq \emptyset$, then 
    $U\subseteq V$.
  \item[(b)] If $U\in \Comp(f^{-n}(V))$ and $U\cap V\neq \emptyset$,
    then for all $w,w'\in U,$ 
    \begin{displaymath}
      |(f^n)'(w)|\ge \lambda
\  \textrm{ and } \
       \frac{|(f^n)'(w)|}{|(f^n)'(w')|}\le K. 
    \end{displaymath}
  \item[(c)] $\overline{B(z,r)}\subset U\subset B(z,\kappa r)\sbt\C\sms\overline{\PS(f)}$.
\end{itemize}
\end{thm}

\fr Each nice set canonically gives rise to a countable alphabet conformal iterated function system in the sense considered in the previous sections of the present paper. Namely, put 
$$
\Comp_*(V)=\bu_{n=1}^\infty\Comp(f^{-n}(V)).
$$
For every $U\in \Comp_*(V)$ let $\tau_V(U)\ge 1$ the unique integer $n\ge 1$
such that $U\in \Comp(f^{-n}(V))$. Put further
$$
\phi_U:=f_U^{-\tau_V(U)}:V\to U
$$
and keep in mind that
$$
\phi_U(V)=U.
$$
Denote by $E_V$ the subset of all elements $U$ of $\Comp_*(V)$ such that 

\sp\begin{itemize}

\item[(a)] $\phi_U(V)\sbt V$, 

\sp\item[(b)] $f^k(U)\cap V=\es$ \, for all \, $k=1,2,\ldots,\tau_V(U)-1$.
\end{itemize}

\sp\fr The collection 
$$
\cS_V:=\{\phi_U:V\to V\}
$$
of all such inverse branches forms obviously a conformal iterated function system in the sense considered in the previous sections of the present paper. In other words, the elements of $\cS_V$ are formed by all holomorphic inverse branches of the first return map $f_V:V\to V$. In particular, $\tau_V(U)$ is the first return time of all points in $U=\phi_U(V)$ to $V$. We define the function $N_V:E_V^\N\to\N_1$ by setting
$$
N_V(\om):=\tau_V(\om_1).
$$
Let 
$$
\pi_V:E_V^\N\to\oc
$$
be the canonical projection induced by the iterated function system $\cS_V$. Let
$$
J_V=\pi_V\(E_V^\N\)
$$
be the limit set of the system $\cS_V$. Clearly
$$
J_V\sbt J(f).
$$
It is immediate from our definitions that
$$
\tau_V(\pi(\om))=N_V(\om)
$$
for all $\om\in E_V^\N$. It is a general fact from abstract ergodic theory that $\mu_{t,V}$, the conditional measure of $\mu_t$ on $V$ is $f_V$--invariant and ergodic. It is clear that $\mu_{t,V}$ is the (only) equilibrium state of the H\"older continuous summable potential
$$
\tilde\psi_{t,V}:=\psi_{t,V}-\P(\psi_t)\tau_V:J_V\to\R,
$$
where 
$$
\psi_{t,V}(x):=\sum_{j=0}^{\tau_V(x)-1}\psi_t\circ f^j(x).
$$
Since the point $z$ is recurrent, $z\in J_V$, and since $z\in V$, the Full Thin Annuli Property of measure $\mu_t$ will follow from Theorem~\ref{t1alr4} provided that condition \eqref{1_2016_02_22} and geometric irreducibility are verified. But the former follows from the assumption that $t>\rho_f/\alpha$ while the latter holds since the Julia set is not contained in any real analytic 
curve.
\end{proof}

\end{document}